%% file: arxiv_v5.tex
\documentclass[12pt]{article}

\usepackage[mathscr]{euscript}
\usepackage[makeroom]{cancel}
\usepackage[font={small,it,singlespacing}]{caption}
\usepackage{amsmath, amssymb, amsthm, amsbsy, amsfonts, bigints, mathtools, cuted}
\usepackage{hyperref, enumerate, verbatim, sectsty, colortbl, fancyhdr}
\usepackage[round, sort, numbers]{natbib}
\setcitestyle{square}
\usepackage{graphicx, float, pdflscape,rotating, pspicture, epsfig, subcaption, threeparttable, multirow, multicol, afterpage,lscape}
\usepackage[nodisplayskipstretch]{setspace}
\usepackage[compact]{titlesec}
\usepackage[noend]{ algpseudocode}
\usepackage{algorithm, enumitem}
\usepackage{appendix}
\definecolor{Gray}{gray}{0.9}
\usepackage[belowskip=1pt,aboveskip=1pt]{caption}
\usepackage[compact]{titlesec}
\usepackage{multicol}

\hypersetup{colorlinks,citecolor=blue}
\graphicspath{{../IEEE-TKDE/}}

\newcommand{\R}{\mathcal{R}}
\newcommand{\T}{\mathcal{T}}
\newcommand{\x}{\mathbf{x}}

\newcommand{\s}{\mathbf{s}}
\newcommand{\e}{\mathbf{e}}

\newcommand{\tthh}{\mbox{\textsuperscript{th}}}

\newcommand{\maxx}{\underset{
\begin{subarray}{c}
  \scriptstyle{\mathbf{x},\mathbf{x}'} \\
  \scriptstyle{d(\mathbf{x},\mathbf{x}')=1}
  \end{subarray}}{\mbox{max}}}

\newcommand{\maxxu}{\underset{
 \begin{subarray}{c}
  \scriptstyle{\mathbf{x},\mathbf{x}',\s^*\in\mathcal{S}} \\
  \scriptstyle{d(\mathbf{x},\mathbf{x}')=1}
  \end{subarray}}{\mbox{max}}}

 \newcommand{\maxxtex}{\mbox{max}_{\mathbf{x},\mathbf{x}',d(\mathbf{x},\mathbf{x}')=1}}
\newcommand{\maxxutex}{\mbox{max}_{\mathbf{x},\mathbf{x}',d(\mathbf{x},\mathbf{x}')=1, \s^*\in\mathcal{S}}}

\theoremstyle{definition}
\newtheorem{thm}{Theorem}
\newtheorem{defn}[thm]{Definition}
\newtheorem{lem}[thm]{Lemma}
\newtheorem{cor}[thm]{Corollary}

\newtheorem{claim}[thm]{Claim}

\begin{document}
\title{\large{\textbf{Generalized Gaussian Mechanism for Differential Privacy}}}
\author{\small{Fang Liu}\footnote{Fang Liu is Associate Professor in the Department of Applied and Computational Mathematics and Statistics, University of Notre Dame, Notre Dame, IN 46556 ($^{\ddag}$E-mail: fang.liu.131@nd.edu). The work is supported by the NSF Grant  1546373 and the University of Notre Dame Faculty Research Support Program Initiation Grant.}}
\date{}

\maketitle

\begin{abstract}
Assessment of disclosure risk is of paramount importance in the research and applications of data privacy techniques. The concept of  differential privacy (DP) formalizes privacy in probabilistic terms and provides a robust concept for privacy protection  without making  assumptions about the background knowledge of adversaries.  Practical applications of DP involve development of DP mechanisms to release results at a pre-specified privacy budget. In this paper, we generalize the widely used Laplace mechanism to  the  family of generalized Gaussian (GG) mechanism based on the $l_p$ global sensitivity  of statistical queries. We explore the theoretical requirement for the GG  mechanism to reach DP at prespecified privacy parameters, and investigate the connections and differences between the GG mechanism and the Exponential mechanism based on the GG distribution We also present a lower bound on the scale parameter of the Gaussian mechanism of $(\epsilon,\delta)$-probabilistic DP as a special case of the GG  mechanism, and compare the  statistical utility of the sanitized results in the tail probability and dispersion in the Gaussian and Laplace mechanisms. Lastly, we apply the GG mechanism in 3 experiments (the mildew, Czech, adult data), and compare the accuracy of sanitized results   via the $l_1$ distance and Kullback-Leibler divergence and examine how sanitization affects the prediction power of a classifier constructed with  the sanitized data in the adult experiment.

\end{abstract}

\section{Introduction}\label{sec:introduction}
When releasing information publicly from a database or sharing data with collaborators, data collectors are always concerned about  exposing sensitive personal information of individuals who contribute to the data. Even with key identifiers removed, data users may still identify a participant in a data set such as via linkage with public information.   Differential privacy (DP) provides a strong privacy guarantee to data release without making  assumptions about the background knowledge or behavior of data users \cite{dwork2006calibrating, dwork2008, dwork2011differential}. For a given privacy budget, information released via a differentially private mechanism guarantees  no additional personal information of an individual in the data  can be inferred, regardless how much background information data users already possess about the individual. DP has spurred a great amount work in the development of differentially private mechanisms to release results and data, including  the Laplace mechanism \cite{dwork2006calibrating},  the Exponential mechanism \cite{mcsherry2007mechanism, mcsherry2009privacy}, the medium mechanism \cite{roth2010median}, the multiplicative weights mechanism \cite{multiplicative}, the geometric mechanism \cite{geometric}, the staircase mechanism \cite{staircase}, the Gaussian mechanism \cite{privacybook}, and applications of DP for private and secure inference in a Bayesian setting \cite{privateBayes}, among others. 

In this paper, we unify the Laplace mechanism and the Gaussian mechanism in the framework of a general  family, referred to as the generalized Gaussian (GG) mechanism. The GG mechanism is based on the $l_p$ global sensitivity (GS) of queries,  a generalization of the $l_1$ GS. We demonstrate the nonexistence of a scale parameter  that would lead to a GG mechanism of pure $\epsilon$-DP in the case of $p\ne1$ if the results to be released are unbounded,  but suggest the GG mechanism of $(\epsilon, \delta)$-probabilistic  DP (pDP) as an alternative in such cases. For bounded data we  introduce the truncated GG mechanism and the boundary inflated truncated GG mechanism that satisfy pure $\epsilon$-DP.   We investigate the connections between the GG mechanism and the Exponential mechanism when the utility function in the latter is based on the Minkowski distance, and establish the relationship between the sensitivity of the utility function in the Exponential mechanism and the  $l_p$ GS of queries. We then take a closer look at the Gaussian mechanism (the GG mechanism of order 2), and derive a lower bound on the scale parameter that delivers $(\epsilon,\delta)$-pDP. The bound is tighter  than the bound to satisfy $(\epsilon,\delta)$-approximate DP (aDP) in the Gaussian mechanism \cite{privacybook}, implying less noise being injected in the sanitized results.  We compare the  utility of sanitized results, in terms of the tail probability and dispersion or mean squared errors (MSE), from independent applications of the Gaussian mechanism and the Laplace mechanism. Finally, we run  3 experiments on the mildew, Czech, and adult data, respectively, and sanitize the count data via the Laplace mechanism, the Gaussian mechanisms of $(\epsilon,\delta)$-pDP and $(\epsilon,\delta)$-aDP. We compare the accuracy of sanitized results in terms of the $l_1$ distance and Kullback-Leibler divergence from the original results, and examine how sanitization affects the prediction accuracy of support vector machines constructed with the sanitized data in the adult experiment.

The rest of the paper is organized as follows. Section \ref{sec:GGM}  defines the $l_p$ GS and presents  the GG mechanism of $(\epsilon,\delta)$-pDP, the truncated GG mechanism, and the boundary inflated truncated GG mechanism that satisfy pure $\epsilon$-DP. It also connects and differentiates between the GG mechanisms and the Exponential mechanism when the utility function in the latter is based the Minkowski distance.  Section \ref{sec:gaussian} take a close look at the Gaussian mechanism of $(\epsilon,\delta)$-pDP, and compares it with the Gaussian mechanism of $(\epsilon,\delta)$-aDP. It also compares the tail probability and the dispersion of the noises injected via the Gaussian mechanism of $(\epsilon,\delta)$-pDP and the Laplace mechanism. Section \ref{sec:experiments} presents the findings from the 3 experiments.  Concluding remarks  are given in Section \ref{sec:discussion}.

\section{Generalized Gaussian Mechanism} \label{sec:GGM}
\subsection{differential privacy (DP)}
DP was proposed and formulated  in  Dwork  \cite{dwork2006} and Dwork et al. \cite{dwork2006calibrating}. A perturbation algorithm $\R$ gives $\epsilon$-differential privacy if for all data sets $(\x,\x')$ that differ  by only one individual ($d(\x,\x')=1$), and all possible query results $Q\subseteq \T$ to query $\s$ ($\T$ denotes the output range of  $\R$),
\begin{equation}\label{eqn:dp}
\left|\log\left(\frac{\Pr(\R( \s(\x)) \in Q)}{\Pr(\R( \s(\x'))\in Q)} \right)\right|\le\epsilon,
\end{equation}
\noindent where $\epsilon>0$ is the privacy $``$budget$"$ parameter. $\s$ refers to queries about data $\x$ and $\x'$, we also use it to denote the query results (unless stated otherwise, the domain of the query results is the set of all real numbers).  $d(\x,\x')=1$ is often defined in two ways in the DP community: $\x$ and $\x'$ are of the same size and differ in exactly one record (row) in at least one attributes (columns); and $\x$ is exactly the same as $\x'$ except that it has one less (more) record. Mathematically,  Eqn (\ref{eqn:dp}) states that the probabilities of  obtaining the same  query result perturbed via $\R$ are roughly the same regardless of whether the query is sent to $\x$ or $\x'$.   In layman's terms, DP implies the chance an individual will be identified based on the perturbed  query result is very low since the query result would be about the  same  with or without the individual in the data. The degree of $``$roughly the same$"$ is determined by the privacy budget $\epsilon$. The lower  $\epsilon$ is,  the more similar the probabilities of obtaining the same query results from $\x$ and $\x'$ are.   DP provides a strong and robust privacy guarantee in the sense that it does not assume anything regarding the background knowledge or the behavior on data users.

In addition to the $``$pure$"$ $\epsilon$-DP in Eqn (\ref{eqn:dp}),  there are softer versions of DP, including the  $(\epsilon,\delta)$-approximate DP (aDP) \cite{dwork2006delta}, the $(\epsilon,\delta)$-probabilistic DP (pDP) \cite{onthemap}, the $(\epsilon,\delta)$-random DP (rDP) \cite{randomDP}, and the $(\epsilon,\tau)$-concentrated DP (cDP) \cite{cPD}.  In all the relaxed versions of DP,  one additional parameter is employed  to characterize the amount of relaxation on top of the privacy budget $\epsilon$. Both the  $(\epsilon,\delta)$-aDP and the $(\epsilon,\delta)$-pDP reduce to $\epsilon$-DP when $\delta=0$, but are different with respect to the interpretation of $\delta$. In $(\epsilon,\delta)$-aDP,
\begin{equation}\label{eqn:adp}
\Pr(\R( \s(\x)) \in Q)\le e^{\epsilon}\Pr(\R( s(\x'))\in Q) + \delta;
\end{equation}
while a perturbation algorithm $\R$ satisfies $(\epsilon,\delta)$-pDP if
\begin{equation}\label{eqn:pdp}
\Pr\left(\left|\log\left(\frac{\Pr(\R( \s(\x)) \in Q)}{\Pr(\R( \s(\x'))\in Q)} \right)\right|>\epsilon\right)\le\delta;
\end{equation}
that is, the probability of $\R$ generating an output belonging to the disclosure set is bounded below $\delta$, where the disclosure set contains all the possible outputs that leak information for a given privacy  budget $\epsilon$. The fact that probabilities are within $[0,1]$ puts constraints on the values of $\epsilon,\Pr(\mathcal{R}(\s(\x')\in Q)$, and $\delta$  in the framework of $(\epsilon,\delta)$-aDP. By contrast, $(\epsilon,\delta)$-pDP seems to be less constrained and more intuitive with its probabilistic flavor. When $\delta$ is small, $(\epsilon,\delta)$-aDP and  $(\epsilon,\delta)$-aDP are roughly the same.   The $(\epsilon,\delta)$-rDP is also a probabilistic relaxation of DP; but it differs from $(\epsilon,\delta)$-pDP in that the probabilistic  relaxation  is with respect to data generation.  In $(\epsilon,\tau)$-cDP, privacy cost is treated as a random variable with an expectation of $\epsilon$ and the probability of the actual cost $>\epsilon$)$>a$  is bounded by $e^{-(a/\tau)^2/2}$. The $(\epsilon,\tau)$-cDP, similar to the $(\epsilon,\delta)$-pDP, relaxes the satisfaction of DP with respect to $\R$ and is broader in scope.

\subsection{ \texorpdfstring{$l_p$}{} global sensitivity} 
\begin{defn}\label{def:Lpgs}
For all $(\x,\x')$ that is $d(\x,\x')=1$, the $l_p$-global sensitivity  (GS) of query $\s$ is
\begin{equation}\label{eqn:Lpgs}
\Delta_p=\maxx\|\s(\x)-\s(\x')\|_p=\left(\textstyle\sum_{k=1}^r\!\left|s_k(\x)-s_k(\x')\right|^p\right)^{1/p}\mbox{for integer } p\!>\!0.
\end{equation}
\end{defn}
\noindent   In layman's term, $\Delta_{p}$  is the maximum difference measured by the Minkowski distance in query results $\s$  between two neighboring data set $\x,\x'$ with  $d(\x,\x')=1$. The sensitivity is $``$global$"$ since it is defined for all possible data sets and all possible ways that $\x$ and $\x'$ differ by one.  The higher $\Delta_{p}$ is, the more disclosure risk there is on the individuals from releasing the original query results $\s$. The $l_p$ GS is a key concept in the construction of the generalized Gaussian mechanism in Section \ref{sec:GGM}.

The $l_p$ GS is a generalization of the $l_1$ GS \cite{dwork2006, dwork2006calibrating} and the $l_2$ GS \cite{privacybook}. The $``$difference$"$ between  $\s(\x)$ and $\s(\x')$ measured by  $\Delta_{1}$ is the largest among all $\Delta_{p}$ for $p\ge1$ since that $\|\s\|_{p+a} \le \|\s\|_p$  for any real-valued vector $\s$ and $a \ge 0$. In addition, $\Delta_{1}$  is  also the most $``$sensitive$"$ measure given that the rate of change with respective to any $s_k$  is the largest among all $p\ge1$.  When $s$ is a scalar, $\Delta_{p}=\Delta_{1}$ for all $p>0$. When $\s$ is multi-dimensional,  an easy upper bound for $l_1$ GS $\Delta_{1}$  is  $\sum_{k=1}^r\Delta_{1,k},$ the sum of the $l_1$ GS of each element $k$ in $\s$, by the triangle inequality. Lemma \ref{lem:GSp1} gives an upper bound on $\Delta_{p}$ for a general $p$ that includes $p=1$ as a special case (the proof is provided in Appendix \ref{app:GSp1}).
\begin{lem}\label{lem:GSp1}
$\left(\sum_{k=1}^r\Delta_{1,k}^p\right)^{1/p}$ is  an upper bound for $\Delta_{p},$ where $\Delta_{1,k}$ is the $l_1$ GS of $s_k$. 
\end{lem}
\noindent The upper bound given in Lemma \ref{lem:GSp1} can be  conservative in cases where the change from  $\x$ to $\x'$ does not necessarily alter every entry in the multidimensional $\s$. For example, the $l_p$ GS  of releasing a histogram with $r$ bins  is  1 (if $d(\x,\x')=1$ is defined as $\x'$ is one record less/more  than $\x$). In other words, the GS is not $r^{1/p}$ even though there are $r$ counts in the released histogram, but is the same as in releasing a single cell because removing one record only alters the count in a single bin.

It is obvious that each element $s_k$ in $\s$ for $k=1,\ldots,r$ needs to be bounded to obtain a finite $\Delta_p$. The most extreme case is the change from  $\x$ to $\x'$ makes  $s_k$ jump from one extreme to the other, implying the range of $s_k$ can be used as an upper bound for $\Delta_{k,1}$, which, combined with Lemma \ref{lem:GSp1}, leads to the following claim.
\begin{claim}\label{cla:bounds vs GS}
Denote the bounds of statistic $s_k$ by $[c_{k0},c_{k1}]$, both of which are finite.  The GS $\Delta_k\le c_{k1}-c_{k0}$ and the GS for $\s=\{s_k\}_{k=1,\ldots,r}$  is $\Delta_p\le \left(\sum_{k=1}^r (c_{k1}-c_{k0})^p\right)^{1/p}$.
\end{claim}

\subsection{generalized Gaussian distribution}
The GG mechanism is defined based on the GG distribution GG$(\mu,b,p)$ with location parameter $\mu$, scale parameter $b>0$, shape parameter $p>0$. The probability density function (pdf)  is
$$f(x|\mu,b,p) = \frac{p}{2b\Gamma(p^{-1})}\exp\left\{\left(\frac{|x-\mu|}{b}\right)^p\right\}.$$
The mean and variance of $x$ are $\mu$ and $b^2\Gamma(3/b)/\Gamma(1/b)$, respectively. ($\Gamma(t)=\int_0^{\infty} x^{t-1}e^{-x}dx$  is the Gamma function). When $p=1$, the GG distribution is the Laplace distribution with mean $\mu$ and variance  $2b^2$;  when $p=2$, the GG distribution becomes the Gaussian distribution with mean $0$ and variance $b^2/2$.
\begin{figure}[!hbt] \begin{center}
\includegraphics[scale=0.6]{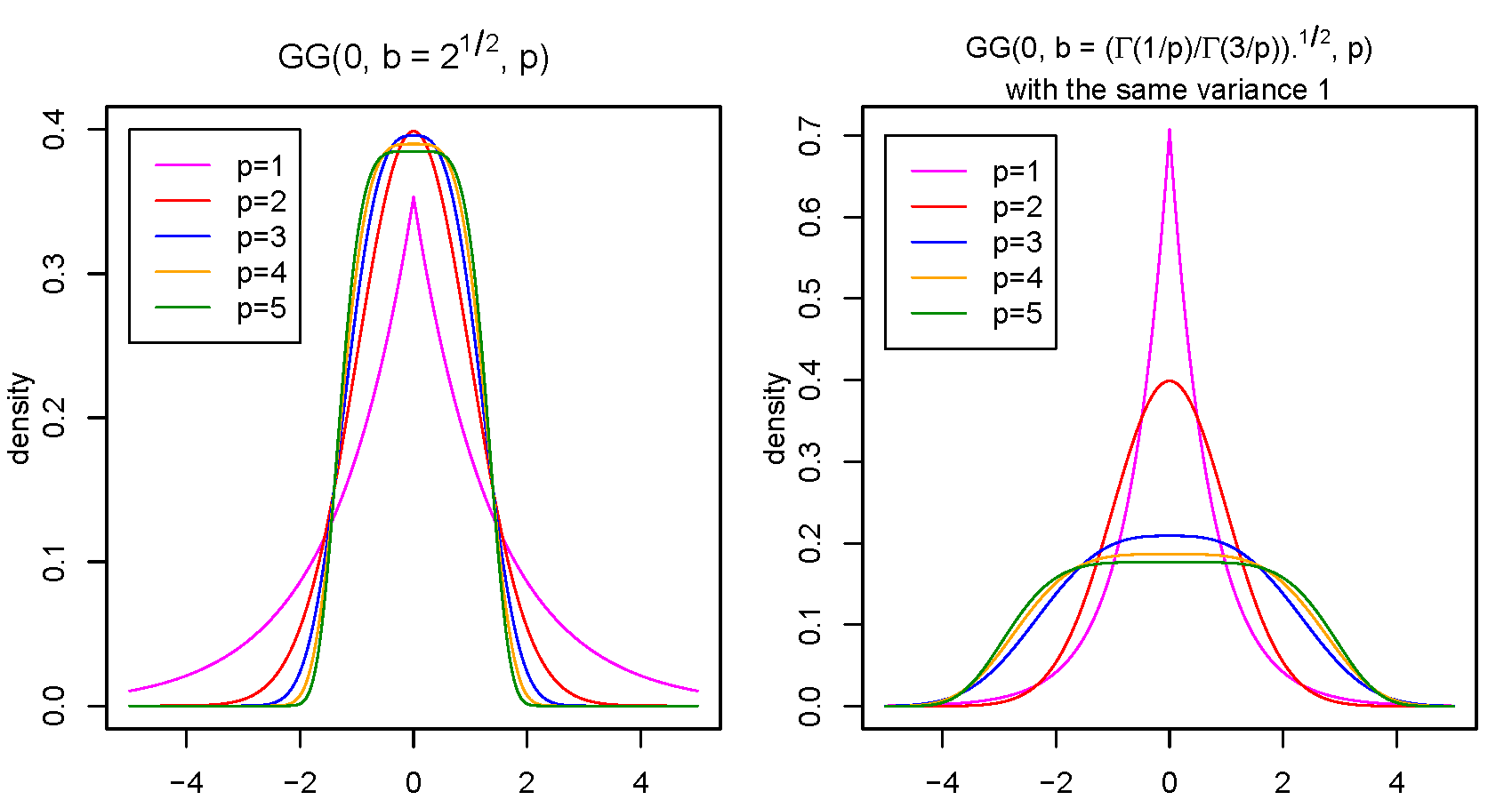}
\caption{Density of GG distributions}\label{fig:GG}
\end{center}\end{figure}
Figure \ref{fig:GG} presents some examples of the GG distributions at different $p$. All the distributions in the left plot have the same scale $b=\sqrt{2}$ and location $0$, and  those in the right plot have the same variance $1$ and location $0$. When the scale parameter is the same (the left plot), the distributions become less spread as $p$ increases, and the Laplace distribution ($p=1$) looks very different from the rest. When the variance is the same (the right plot), the Laplace distribution is the most likely to generate values that are close to the mean, followed by the Gaussian distribution ($p=2$).

\subsection{GG mechanism  of \texorpdfstring{$\epsilon$}{}-DP}\label{sec:GGMDP}
We first examine the GG mechanism of $\epsilon$-DP with the domain for $s_k^*$ defined on $(-\infty,\infty)$ for $k=1,\ldots,r$.   $\s$ needs to bounded to calculate the $l_p$ GS, but the bounding requirement does not necessarily goes into  formulating the GG distribution for the GG mechanism in the first place. If bounding for $\s^*$ is necessary, it can be incorporated in a post-hoc manner after being generated from the GG mechanism. A well-known example is the Laplace mechanism. It employs a Laplace distribution defined on $(-\infty,\infty)$, though its scale parameter $b=\Delta_1/\epsilon$ requires $\s$ to be bounded for $\Delta_1$ to be calculated.

Eqn (\ref{eqn:GGD}) presents the GG distribution from which sanitized  $\s^*$ would be generated to satisfy $\epsilon$-DP,  assuming $b$ exists.
\begin{align}
f(\s^*)&\propto e^{\left(\|\s^*-\s\|_p/b\right)^p}\propto\textstyle\prod_{k=1}^r \exp\{-(|s^*_k-s_k|/b)^p\}\notag\\
&=\textstyle\prod_{k=1}^r \frac{p}{2b\Gamma(p^{-1})}\exp\{(|s^*_k-s_k|/b)^p\} = \textstyle\prod_{k=1}^r \mbox{GG}(s_k,b,p) \label{eqn:GGD}
\end{align}
\begin{claim}\label{cla:lowerb}
There does not exist a lower bound on $b$ for the GG distribution in Eqn (\ref{eqn:GGD}) when $p\ne1$ that generates  $\s^*$ with $\epsilon$-DP.  When $p=1$, the lower bound on $b$ that leads to $\epsilon$-DP is $\epsilon^{-1}\Delta_1$.
\end{claim}
Appendix \ref{app:GGM} lists the detailed steps that lead to Claim \ref{cla:lowerb}. In brief, to achieve $\epsilon$-DP, we need
$b^{-p} \left(\textstyle\sum_{k=1}^r\sum_{j=1}^{p-1}\!(_j^p)|s^*_k-s_k|^{p-j}\Delta_{1,k}^j+\Delta_{p}^p\right)\le\epsilon$ (Eqn \ref{eqn:fourth}). However,  this inequality depends on the random GG noise $e_k=s^*_k-s_k$ for $k=1,\ldots,r$, the support of which is $(-\infty,\infty)^r$. In other words, there does not exist a random noise-free solution on $b$, unless $p=1$ in which case the inequality no longer involves the error terms and the GG mechanism reduces to the familiar Laplace mechanism of $\epsilon$-DP.   We propose two approaches to fix the problem and achieve DP through the GG mechanism. The first approach leverages the bounding requirement  for $\s$ and builds in the requirement in the GG distribution in the first place to generate $\s^*$ with $\epsilon$-DP, assuming that $\s^*$ and $\s$ share the same bounded domain (Section  \ref{sec:boundedGGM}). The second approach still uses the GG distribution in Eqn (\ref{eqn:GGD}) to sanitize $\s$, only satisfying $(\epsilon,\delta)$-pDP instead of the pure $\epsilon$-DP (Section \ref{sec:GGMpDP}).  The sanitized $\s^*$ can be bounded in a post-hoc manner, as needed.

\subsection{truncated GG mechanism and boundary inflated truncated GG mechanism of \texorpdfstring{$\epsilon$}{}-DP}\label{sec:boundedGGM}
\begin{defn} \label{def:tGGM}
Denote the bounds on query result $\s$ by $[c_{k0},c_{k1}]_{k=1,\ldots,r}$. For integer $p\ge1$, the  truncated GG mechanism of order $p$ generates $\s^*\!\!\in\!\![c_{k0},c_{k1}]_{k=1,\ldots,r}$  with $\epsilon$-DP by drawing  from the truncated GG distribution
\begin{align}
&f(\s^*|c_{k0}\!\le\!s^*_k\! \le\! c_{k1},\forall\; k=1,\ldots,r)=\prod_{k=1}^r\frac{p\exp\{(|s^*_k-s_k|/b)^p\}}{2b\Gamma(p^{-1}) A(s_k, b, p)} \mbox{ with scale parameter }\label{eqn:tGGD}\\
&\qquad\qquad b \ge \left(2\epsilon^{-1}\!\!\left(\!\sum_{k=1}^r\sum_{j=1}^{p-1}(_j^p)|c_{k1}-c_{k0}|^{p-j}\Delta_{1,k}^j+\Delta_{p}^p\right)\!\right)^{1/p}\!\!\!\!\!\!\!,
\label{eqn:b1}
\end{align}
where $A(s_k, b,p)\!=\!\Pr(c_{k0}\!\le\!s^*_k\!\le\! c_{k1};  s_k, b,p)\!=\!(\Gamma(p^{-1}))^{-1}(\gamma[p^{-1},(c_{k1}-s_k)/b] +\gamma[p^{-1},(s_k-c_{k0})/b])$ ($\gamma$ is the lower incomplete gamma function), $\Delta_{1,k}$ is the $l_1$ GS of $s_k$, and $\Delta_{p}$  is the $l_p$ GS  of $\s$.
\end{defn}
The proof of $\epsilon$-DP of the truncated GG mechanism is given in Appendix  \ref{app:tGGM}. The truncated GG mechanism perturbs each element in  $\s$ independently; thus  Eqn (\ref{eqn:tGGD}) involves the product of $r$ independent density functions.   Though the closed interval $[c_{k0},c_{k1}]$ is used to denote the bounds on $s_k$, Definition \ref{def:tGGM} remains the same regardless of whether the interval is closed, open, or half-closed since the GG distribution is defined on a continuous domain. If $s_k$ is discrete in nature such as counts, post-hoc rounding on perturbed  $\s_k^*$ can be applied. The lower bound on $b$ in Eqn  (\ref{eqn:b1}) depends on $\Delta_{p}$. We may apply Lemma \ref{lem:GSp1} and set $\Delta_{p}^p$ at its upper bound $\sum_{k=1}^r\Delta^p_{1,k}$ to obtain a less tight bound on $b$.
\begin{align}\label{eqn:b2}
\!\!\!\!\!b\! \ge\!\left(2\epsilon^{-1}\textstyle\!\left(\!\sum_{k=1}^r\sum_{j=1}^p(_j^p)|c_{k1}-c_{k0}|^{p-j}\Delta_{1,k}^j\right)\!\right)^{1/p}\!\!\!.
\end{align}

\begin{defn} \label{def:bitGGM}
Denote the bounds on query result $s_k$ by $[c_{k0},c_{k1}]$ for $k=1,\ldots,r$. For integer $p\ge1$, the $p\tthh$ order  boundary inflated truncated (BIT) GG mechanism sanitizes $\s$ with $\epsilon$-DP by drawing perturbed  $\s^*$ from the following piecewise distribution
\begin{align}\label{eqn:thGGD}
f(\s^*|c_{k0}\!\le\!s^*_k\! \le\! c_{k1},\forall\; k=1,\!\ldots,\!r)\!=\! \textstyle\prod_{k=1}^r\!\left\{\!p_k^{\mathrm{I}(s_k^*=c_{k0})}q_k^{\mathrm{I}(s_k^*=c_{k1})} \!\left(\frac{p\exp\{(|s^*_k-s_k|/b)^p\}}{2b\Gamma(p^{-1})}\right)^{\!\!\mathrm{I}(c_{k0}<s_k^*<c_{k1})}\right\}\!,\!\!
\end{align}
where $p_k\!=\!\Pr(s_k^*\!<\! c_{k0}; s_k,p,b)\!=\!\frac{1}{2}\!-\!\gamma(p^{-1}, ((s_k\!-\!c_{k0})/b)^p)(2\Gamma(p^{-1}))^{-1}$ and $q_k=\Pr(s_k^*> c_{k1}; s_k,p,b) =\frac{1}{2}-\gamma(p^{-1}, ((c_{k1}-s_k)/b^p))(2\Gamma(p^{-1}))^{-1}$, $\gamma$ is the lower incomplete gamma function, and $\Gamma$ is the gamma function; and $\mathrm{I}()$ is the indicator function that equals 1 if the argument in the parentheses is true, 0 otherwise.
\end{defn}
In brief, the BIT GG distribution replaces out-of-bound values with the boundary values and keeps the within-bound values as is, leading to  a piecewise distribution. This is in contrast to the truncated GG distribution which  throws away out-of-bound values.  The challenge with perturbing $\s$ directly via Eqn (\ref{eqn:thGGD}) lies in solving for a lower bound $b$ that satisfies $\epsilon$-DP  from
 \begin{equation}\label{eqn:bit.ineq1}
 \log\left|\frac{f(\s^*|c_{k0}\!\le\!s^*_k\! \le\! c_{k1},\forall\; k=1,\ldots,r)}{f(\s^{'*}|c_{k0}\!\le\!s^*_k\! \le\! c_{k1},\forall\; k=1,\ldots,r)}\right|\le\epsilon
\end{equation}
where $\s^*=\{s^*_k\}$ and $\s'^{*}=\{s'^{*}_k\}$ are the sanitized results from data $\x$ and $\x'$ that are $d(\x,\x')=1$, respectively. The lower bound given in Eqns (\ref{eqn:b1}) and \ref{eqn:b2}  can be used when the output subset $Q$ is a subset of $(c_{10}, c_{11})\times\cdots\times(c_{r0}, c_{r1})$ (open intervals). However, when $Q$ is $\{s_k=c_{k0}\;\forall\; k=1,\ldots,r\}$ and  $\{s_k=c_{k1}\;\forall\; k=1,\ldots,r\}$, respectively, there  are no  analytical solutions  on $b$ in either Eqns (\ref{eqn:bit.ineq2}) or  (\ref{eqn:bit.ineq3})
\begin{align}
\!\!\!\!\!&\log\!\left|\textstyle\prod_{i=1}^r\!\frac{1/2\!-\!\gamma(p^{-1}, (\!(s_k-c_{k0})/b)^p)(2\Gamma(p^{-1}))^{-1}}{1/2-\gamma(p^{-1}, ((s'_k-c_{k0})/b)^p)(2\Gamma(p^{-1}))^{-1}}\!\right|\!\le\!\epsilon\label{eqn:bit.ineq2}\\
\!\!\!\!\!&\log\!\left|\textstyle\prod_{i=1}^r\!\frac{1/2\!-\!\gamma(p^{-1}, (\!(s_k-c_{k0})/b)^p)(2\Gamma(p^{-1}))^{-1}}{1/2-\gamma(p^{-1}, ((s'_k-c_{k0})/b)^p)(2\Gamma(p^{-1}))^{-1}}\!\right|\!\le\!\epsilon.\label{eqn:bit.ineq3}
\end{align}
 The most challenging situation is when $Q$  is a mixture set of $(c_{k0}, c_{k1})$, $c_{k0}$, and $c_{k1}$ for different $k=1,\ldots,r$. In summary,  the BIT GG mechanism is not very appealing from a practical standpoint.

\subsection{GG mechanism of  \texorpdfstring{$(\epsilon, \delta)$}{}-pDP}\label{sec:GGMpDP}
The  second approach to obtain a lower bound on the scale parameter $b$ for the GG distribution in Eqn (\ref{eqn:GGD}) when $p\ge 2$ is to employ a soft version of DP. Corollary \ref{cor:GGMpDP} presents a solution on  $b$ that satisfies $(\epsilon,\delta)$-pDP.
\begin{cor} \label{cor:GGMpDP}
If the scale parameter $b$ in the GG distribution in Eqn (\ref{eqn:GGD})  satisfies
\begin{align}
\!\!&\Pr\!\left(\!\textstyle\sum_{k=1}^r\!\sum_{j=1}^{p-1}(_j^p)|s^*_k\!-\!s_k|^{p-j}\Delta_{1,k}^j\! >\!b^p\epsilon\!-\!\Delta_{p}^p\!\right)\!<\! \delta,\label{eqn:GGMpDP}
\end{align}
then the GG mechanism satisfies $(\epsilon,\delta)$-pDP when $p\ge2$.
\end{cor} 
\noindent The proof is straightforward. Specifically, rather than setting the left side of Eqn (\ref{eqn:fourth}) $\le\epsilon$ (i.e. with 100\%), we attach a probability of achieving the inequality, that is, Pr(Eqn (\ref{eqn:fourth})$<\epsilon)>1-\delta$, leading to Eqn (\ref{eqn:GGMpDP}). The $(\epsilon,\delta)$-pDP does not apply to the Laplace mechanism ($p=1$) at least in the framework laid out in Corollary \ref{cor:GGMpDP}. When $p\!=\!1$,  Eqn (\ref{eqn:first}) becomes $b^{-1}\sum_{k=1}^r\!\big||e_k|\!-\!|e_k+d_k|\big|\!\le\! b^{-1}\sum_{k=1}^r|d_k|\!\le\! b^{-1}\Delta_{1}$, which does not involve the random variable $\s^*$; in other words, as long as  $b^{-1}\Delta_{\s,1}\le\epsilon$, the pure $\epsilon$-DP is guaranteed.

Corollary \ref{cor:GGMpDP} does not list a closed-form solution on $b$ as it is likely that  only numerical solutions  exist in most cases. Given that $s^*_k$ is independent across $k=1,\ldots,r$, $a_k= \textstyle\sum_{j=1}^{p-1}(_j^p)|s^*_k-s_k|^{p-j}\Delta_{1,k}^j $ a function of $s^*_k$, is also independent across $k$. Therefore, the problem becomes searching for a lower bound on $b$ where the probability of a sum of $r$ independent variables ($a_1,\ldots,a_r$) exceeding $b^p-\Delta_p^p\epsilon$ is smaller than $\delta$.  If there exists a closed-form distribution function for $\sum_{k=1}^r a_k$,  an exact solution on $b$ can be obtained. When $p=2$, an analytical lower bound $b$ can be obtained (see Section \ref{sec:gaussian}); when $p> 2$  we only manage to obtain the distribution function for $(_j^p)|s^*_k-s_k|^{p-j}\Delta_{1,k}^j$, but not  for $a_k$ or $\sum_{k=1}^r a_k$ at the current stage.  A relatively simple case is when the elements of statistics $\s$ are calculated on disjoint subsets of the original data, thus removing one individual from the data only affects one element out of $r$, $\Delta_1=\Delta_p= \Delta_{1,k'}$, leading to the Corollary \ref{cor:disjoint}.
\begin{cor} \label{cor:disjoint}
When all $r$ elements in $\s$ are based disjoint subsets of the data, the lower bound on $b$ satisfies
$\Pr(\sum_{j=1}^p(_j^p)|s^*_{k'}\!-\!s_{k'}|^{p-j}\Delta_{1,k'}\!>\!b^p\epsilon)\!<\! \delta$, where $k'=\mbox{argmax}_k\Delta_{1,k}$.
\end{cor}
When the query is a histogram, $\Delta_1=\Delta_p= \Delta_{1,k'}=1$, and the lower bound $b$ for $(\epsilon,\delta)$-pDP can be derived from $\Pr(\sum_{j=1}^p(_j^p)|e_{k'}|^{p-j}\!>\!b^p\epsilon)\!<\! \delta$. The proof of \ref{cor:disjoint} is trivial.  With disjoint queries, only one element in $\s$ is affected by changing from $\x$ to $\x'$ while the other $r-1$ elements in Eqn (\ref{eqn:second}) in Appendix \ref{app:GGM} are 0 as $s_k(\x)=s_k(\x')$, and Eqn (\ref{eqn:second}) $=b^{-p}\sum_{j=1}^p(_j^p)|e_{k'}|^{p-j}|d_{k'}|^j\le|b^{-p}\sum_{j=1}^p(_j^p)|e_{k'}|^{p-j}\Delta_{1,k'}$.

\begin{figure}[htb] \begin{center}
\includegraphics[scale=0.6]{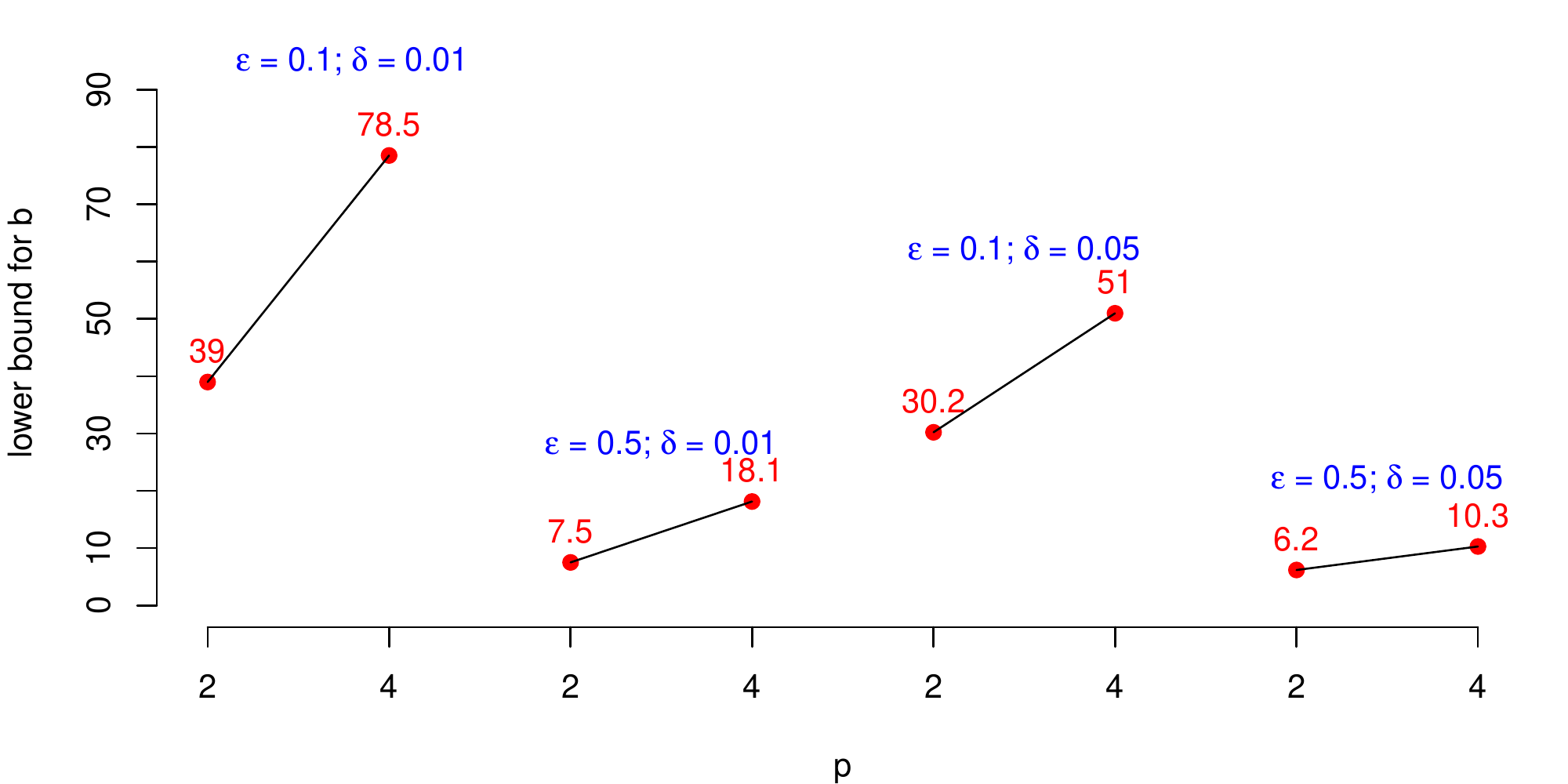}
\caption{Numerical Lower bound on $b$ from Corollary \ref{cor:GGMpDP}}\label{fig:lowerbound}
\end{center}\end{figure}

Numerical approaches can be applied to obtain a lower bound on $b$ when the closed-form solutions are difficult to attain. Figure \ref{fig:lowerbound} depicts the lower bounds on $b$ at different $p$ and $(\epsilon,\delta)$ obtained via the Monte Carlo approach. We set $\Delta_{1,k}$ at $1,0.1,0.05$ for $k=1,2,3$, respectively and applied Lemma \ref{lem:GSp1} to obtain an upper bound on $\Delta_{p}$ for a given $p$ value. As expected, the lower bound on $b$ increases with decreased $\epsilon$ (lower privacy budget) and decreased $\delta$ (reduced chance of failing the pure $\epsilon$-DP).   The results also suggest $b$ increases with $p$ to maintain $(\epsilon,\delta)$-pDP in the examined scenarios.

$\s^*$ sampled from the GG mechanism  of $(\epsilon,\delta)$-pPD in Eqn (\ref{eqn:GGD}) once $b$ is determined -- analytically or numerically -- ranges $(-\infty,\infty)$. To bound $\s^*$, it is straightforward to apply a post processing procedure such as the truncation and the boundary inflated truncation (BIT) procedure \cite{liu2016b}. The truncation procedure throws away the out-of-bounds values and only keeps those in bounds while the BIT procedure sets the out-of-bounds values at the bounds.  If the bounds are noninformative in the sense that the bounds are global and do not contain any data-specific information,  then neither one of the two post-hoc bounding procedures will leak the original information  or compromise the  established $(\epsilon,\delta)$-pDP.

\subsection{Connection between GG mechanism and Exponential Mechanism}\label{sec:connection}
The exponential mechanism  was  introduced by  McSherry and Talwar \cite{mcsherry2007mechanism}. We paraphrase the original definition as follows, covering both discrete and continuous outcomes. Let $\mathcal{S}$ denote the set containing all possible output $\s^{\ast}$. The exponential mechanism releases $\s^{\ast}$  with probability
\begin{equation}\label{eqn:exp}
f(\s^{\ast})=\exp\left(u(\s^{\ast}|\x)\frac{\epsilon}{2\Delta_u}\right)(A(\x))^{-1}
\end{equation}
\noindent to ensure $\epsilon$-DP. $A(\x)$ is a normalizing constant so that $f(\s^{\ast})$ sums or integrates to 1,  and equals to $\sum_{\s^{\ast}\in\mathcal{S}} \exp\left(\!u(\s^{\ast}|\x)\frac{\epsilon}{2\Delta_u}\!\right)$ or $\int_{\s^{\ast}\in\mathcal{S}} \exp\left(\!u(\s^{\ast}|\x)\frac{\epsilon}{2\Delta_u}\!\right)\!d\s^{\ast}$, depending on whether $\mathcal{S}$ is a countable/discrete sample space, or a continuous set, respectively. $u$ is the utility function and assigns a $``$utility$"$ score to each possible outcome $\s^*$ conditional on the original data $\x$, and $\Delta_u=\maxxutex|u(\s^{\ast}|\x)-u(\s^{\ast}|\x')|$ is the maximum change in the utility score across all possible output $\s^*$ and all possible data sets $\x$ and $\x'$ that is $d(\x,\x')=1$.
From a practical perspective, the scores should properly reflect the $``$usefulness$"$ of $\s^*$. For example, $``$usefulness$"$ can be measured the similarity between perturbed  $\s^{\ast}$ and original $\s$ if $\s$ is numerical. The closer $\s^{\ast}$ is to the original $\s$, the larger $u(\s^{\ast}|\x)$ is, and the higher the probability $\s^{\ast}$ will be released.
The Exponential mechanism can be conservative (See Appendix \ref{app:exp}), in the sense that the actual privacy cost is lower than the nominal privacy budget $\epsilon$, or more than necessary amount of perturbation is injected to preserve $\epsilon$-DP. Despite the conservativeness, the Exponential mechanism is a widely used mechanism in DP with its generality  and flexibility as long as the utility function $u$ is properly designed.

When $u$ is defined as the negative $p\textsuperscript{th}$ power of the $p\tthh$-order Minkowski distance between $\s^{\ast}$ and $\s$, that is, $u(\s^{\ast}|\s)\!=\!-\|\s^{\ast}-\s\|_p^p$, the Exponential mechanism generates perturbed  $\s^*$ from the GG distribution
\begin{align}\label{eqn:explp0}
&\textstyle f(\s^{\ast}|\s)\!=\!(A(\s))^{-1}\!\exp\!\left(\!\!-\|\s^{\ast}-\s\|_p^p\frac{\epsilon}{2\Delta_u}\!\right)
=\textstyle(A(\s))^{-1}\!\prod_{k=1}^r\exp\!\left(\!-\frac{|s_k^{\ast}-s_k|^p}{2\Delta_u\epsilon^{-1}}\!\right)\!=\!\prod_{k=1}^r\mbox{GG}(s_k,b,p)
\end{align}
with  $A(\s)\!=\!\left(p^{-1}2b\Gamma(p^{-1})\right)^{\!r}$ and $b^p\!=\!2\Delta_u\epsilon^{-1}$. The scale parameter $b$  in Eqn (\ref{eqn:explp0}) is a function of  the GS  of the utility function $\Delta_u$ and  the privacy budget $\epsilon$. For bounded data $s_k^*\in[c_{k0},c_{k}]$ for $k=1,\ldots,r$, the Exponential mechanism based on the GG distribution is
\begin{align}\label{eqn:explp1}
&f(\s^{\ast}|\s^*\in[\mathbf{c}_0,\mathbf{c}_1])=
(A(\s))^{-1}\textstyle\prod_{k=1}^r(B(s_k))^{-1}\exp\!\left(\!-\frac{|s_k^{\ast}-s_k|^p}{2\Delta_u\epsilon^{-1}}\!\right),
\end{align}
where $B(s_k)=\Pr(s_k^*\in[c_{k0},c_{k}])$ is  calculated from the pdf $\mbox{GG}(s_k,b,p)$. Compared to the truncated GG mechanism in Definition \ref{def:tGGM}, the only difference in the Exponential mechanism in Eqn (\ref{eqn:explp1}) is how the scale parameter $b$ is defined. In  Definition \ref{def:tGGM}, $b$ depends on the GS  of $\s$ ($\Delta_p$) while it is a function of the GS  of the utility function $u$ ($\Delta_u$) in the Exponential mechanism.  Specifically, $b^p\ge2\epsilon^{-1}\Delta_u$ in the Exponential mechanism, and the lower bound on $b$ is given in Eqn (\ref{eqn:b1}) in the GG mechanism. While both mechanisms will lead to the satisfaction of $\epsilon$-DP, the one with a smaller $b$ is preferable at the same $\epsilon$.  The magnitude of  $b$  in each case depends on the bounds of $\s$, and the order $p$, in addition to $\Delta_u$ or $\Delta_{p}$. Though not a direct comparison on $b$, Lemma \ref{lem:deltau.s.relationship} explores the relationship between $\Delta_u$ and $\Delta_{p}$, with the hope to shed light on the comparison of $b$ (the proof is in Appendix \ref{app:deltau.s.relationship}).
\begin{lem}\label{lem:deltau.s.relationship}
Let $[c_{k0},c_{k1}]$ denote the bounds on $s_k$ for $k=1,\ldots,r$.
\begin{enumerate}
\item[a)] When $u=-\|\s^{\ast}-\s\|_1$,  $\Delta_u\le\Delta_{1}$.  Both the GG mechanism and the GG-distribution based Exponential mechanism reduce to the truncated Laplace mechanism with the same $b$.
\item[b)] When $u=-\|\s^{\ast}-\s\|_2^2$, $\Delta_u\le2\sum_{k=1}^r\Delta_{1,k}|c_{k1}-c_{k0}|$.
\item[c)] When $u\!=-\|\s^{\ast}-\s\|_p^p$ for $p\ge3$,  $\Delta_u\!\le\!\sum_{k=1}^r\!\sum_{j=1}^p\! (^p_j)\!\left(\mbox{max}\{|c_{k0}|, |c_{k1}|\}\right)^{p-j}\!\Delta_{1,k}^{(j)}$, where $\Delta^{(j)}_{1,k}=\!\maxxtex|(s_k(\x))^j-(s_k(\x'))^j|$ is $l_1$ GS of $(s_k)^j$.
\end{enumerate}
\end{lem}
As a final note on the GG-distribution based Exponential mechanism, we did not use the negative Minkowski distance directly as the utility function due to a couple of potential practical difficulties with this approach.  First, $\Delta_u$ can be difficulty to obtain. Second, $f(\s^*)\!\propto\!\textstyle\exp\{-\left(\sum_{k=1}^r |s_k^*-s_k|^p\right)^{1/p}\epsilon(2\Delta_u)^{-1}\}$,  does not appear to be associated with any known distributions (except when $p=1$), and additional efforts are required to study the properties of $f(\s^*)$ and to develop an efficient algorithm to draw samples from it.

\section{Gaussian Mechanism}\label{sec:gaussian}
A special case of the GG mechanism is the Gaussian mechanism when $p=2$ that draws $s_k^*$  independently from  a Gaussian distribution with mean $s_k$ and variance $\sigma^2=b^2/2$ for $k=1,\ldots,r$. Applying Eqn (\ref{eqn:tGGD}) with $b$ defined in Eqns (\ref{eqn:b1}) and (\ref{eqn:b2}), we can obtain the truncated Gaussian mechanism of $\epsilon$-DP for bounded $\s\in[c_{10},c_{11}]\times\cdots\times[c_{r0},c_{r1}]$
\begin{align}\label{eqn:GaussianDP}
f(\s^*|\s) &=\textstyle \prod_{k=1}^r\left\{\left(\Phi(c_{k1};\mu, \sigma^2)-\Phi(c_{k0};\mu, \sigma^2)\right)^{-1}\phi(s_k^*;\mu=s_k, \sigma^2=b^2/2)\right\},\mbox{ where}\\
b^2 &\ge 2\epsilon^{-1}\textstyle\left(2\sum_{k=1}^r|c_{k1}-c_{k0}|\Delta_{1,k}+\Delta_{2}^2\right)\ge 2\epsilon^{-1}\textstyle\sum_{k=1}^r\left(2|c_{k1}-c_{k0}|\Delta_{1,k}+\Delta_{1,k}^2\right),\notag
\end{align}
\noindent where $\phi$ and  $\Phi$ are the pdf and the CDF of the Gaussian distribution, respectively.

\begin{figure}[htb]\begin{center}
\includegraphics[scale=0.85]{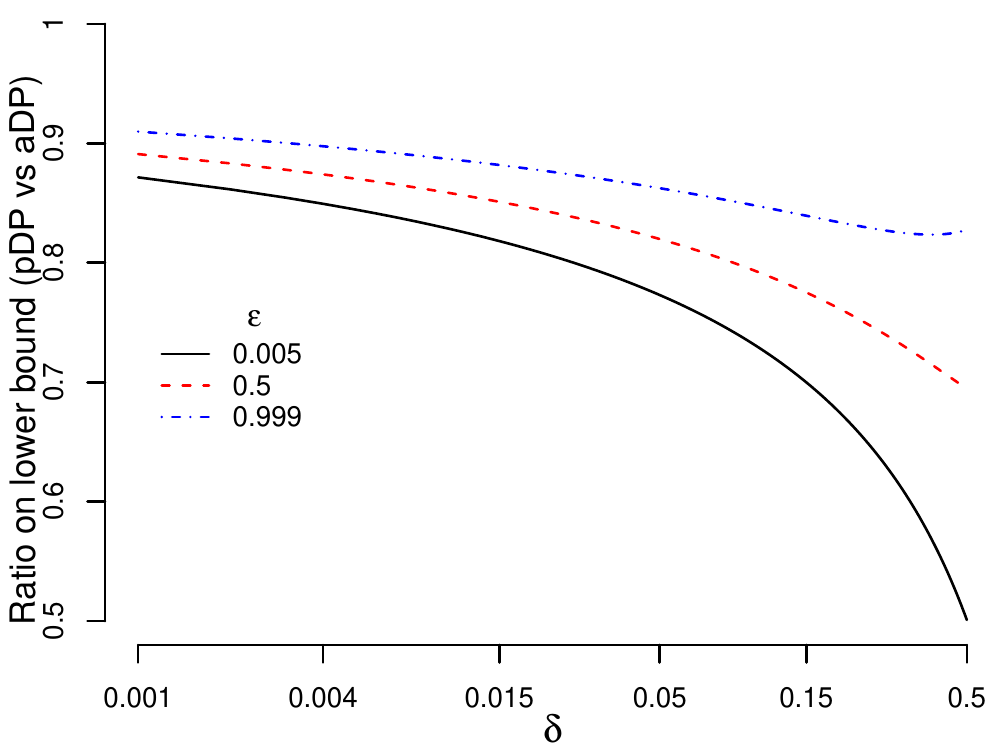}
\caption{Comparison of pDP lower bound (Eqn \ref{eqn:lowerbound2}) vs.  aDP bound (Eqn \ref{eqn:privacybook}) on $\sigma$ in the Gaussian mechanism for $\epsilon<1$  (the aDP bound requires $\epsilon<1$)}\label{fig:lowerbound2}
\end{center}\end{figure}

An analytical solution on the lower bound of $b$ for the Gaussian mechanism of $(\epsilon,\delta)$-pDP is provided in Lemma \ref{lem:lowerbound2} (the proof is provided in Appendix \ref{app:lowerbound2}).
\begin{lem}\label{lem:lowerbound2}
The lower bound on the scale parameter $b$ from the Gaussian mechanism of $(\epsilon,\delta)$-pDP  is
$b\ge2^{-1/2}\epsilon^{-1}\!\Delta_2\!\left(\!\sqrt{(\Phi^{-1}(\delta/2))^2+2\epsilon}-\Phi^{-1}(\delta/2)\!\right)$.
\end{lem}
\noindent Given the relationship between $b$ and the standard deviation  of the Gaussian distribution $\sigma=b/\sqrt{2}$, the lower bound can also be expressed in $\sigma$,
\begin{eqnarray}\label{eqn:lowerbound2}
\sigma \ge  (2\epsilon)^{-1}\Delta_2 \left(\sqrt{(\Phi^{-1}(\delta/2))^2+2\epsilon}-\Phi^{-1}(\delta/2)\right).
\end{eqnarray}
The pDP lower bound given in Eqn (\ref{eqn:lowerbound2}) is different from the lower bound
\begin{equation}\label{eqn:privacybook}
\!\!\!\!\sigma\!>\!\epsilon^{-1}\Delta_2c, \mbox{ with $\epsilon\!\in\!(0,1)$ and $c^2\! >\! 2 \ln(1.25/\delta)$}.
\end{equation}
in Dwork and Roth  \cite{privacybook} for $(\epsilon,\delta)$-aDP (Eqn (\ref{eqn:adp})). The pDP bound in Eqn (\ref{eqn:lowerbound2}) is tighter than the aDP bound in Eqn (\ref{eqn:privacybook}) for the same set of $(\epsilon,\delta)$ (note the interpretation of $\delta$ in pDP and aDP is different, but the DP guarantee is roughly the same when $\delta$ is small). In addition, the pDP bound  does not constrain $\epsilon$ to be $<1$ as required in the aDP bound. Figure \ref{fig:lowerbound2} compares the two two lower bounds at several $\epsilon\in (0,1)$ and $\delta\in(0,0.5)$. As observed, the ratio between the aPD vs. pDP lower bounds  is always $< 1$ for the same $(\epsilon,\delta)$. The smaller $\epsilon$ is, or the larger $\delta$ is, the smaller the ratio is and the larger the difference is between the two bounds.

Dwork and  Roth \cite{privacybook} list several advantages of the Gaussian noises, such as the Gaussian noise is a $``$familiar$"$ type of noise as many noise sources in real life can be well approximated by Gaussian distributions; the sum of Gaussian variable is still a  Gaussian; and finally, in the case of multiple queries or when $\delta$ is small, the pure-DP guarantee in the Laplace mechanism  and the pDP guarantee in the Gaussian mechanism see minimal difference. A theoretical disadvantage to Gaussian noise is that it does not guarantee DP in some cases (e.g.,  Report Noisy Max)\cite{privacybook}.

\begin{figure}[htb] \begin{center}
\includegraphics[scale=0.75]{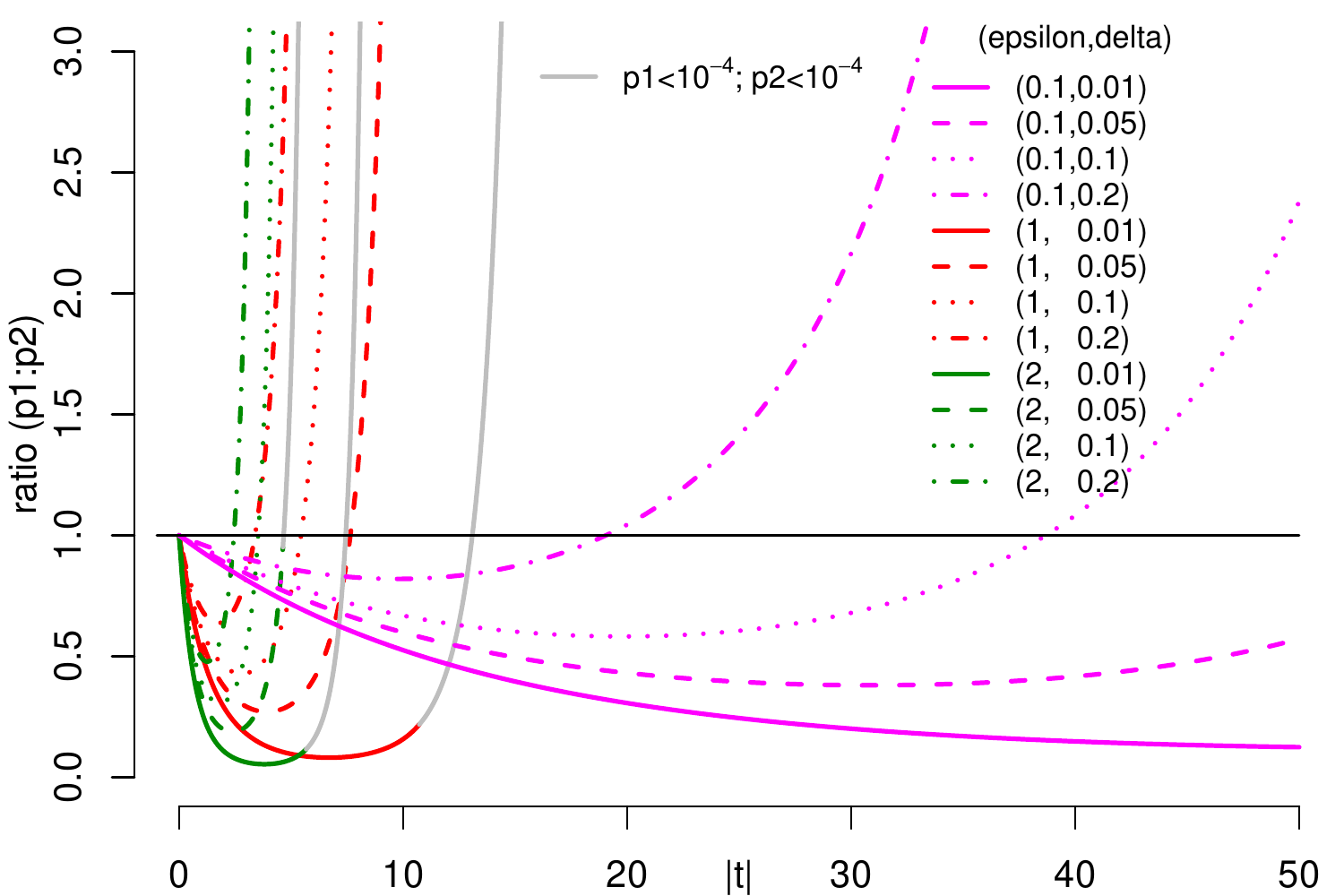}
\caption{Ratio on the tail probabilities $p_1:p_2$ (the gray curves represent the unlikely cases where both $p_1$ and $p_2$ are $<10^{-4}$)}  \label{fig:tail}
\end{center}\end{figure}

\begin{figure}[htb] \begin{center}
\includegraphics[scale=0.8]{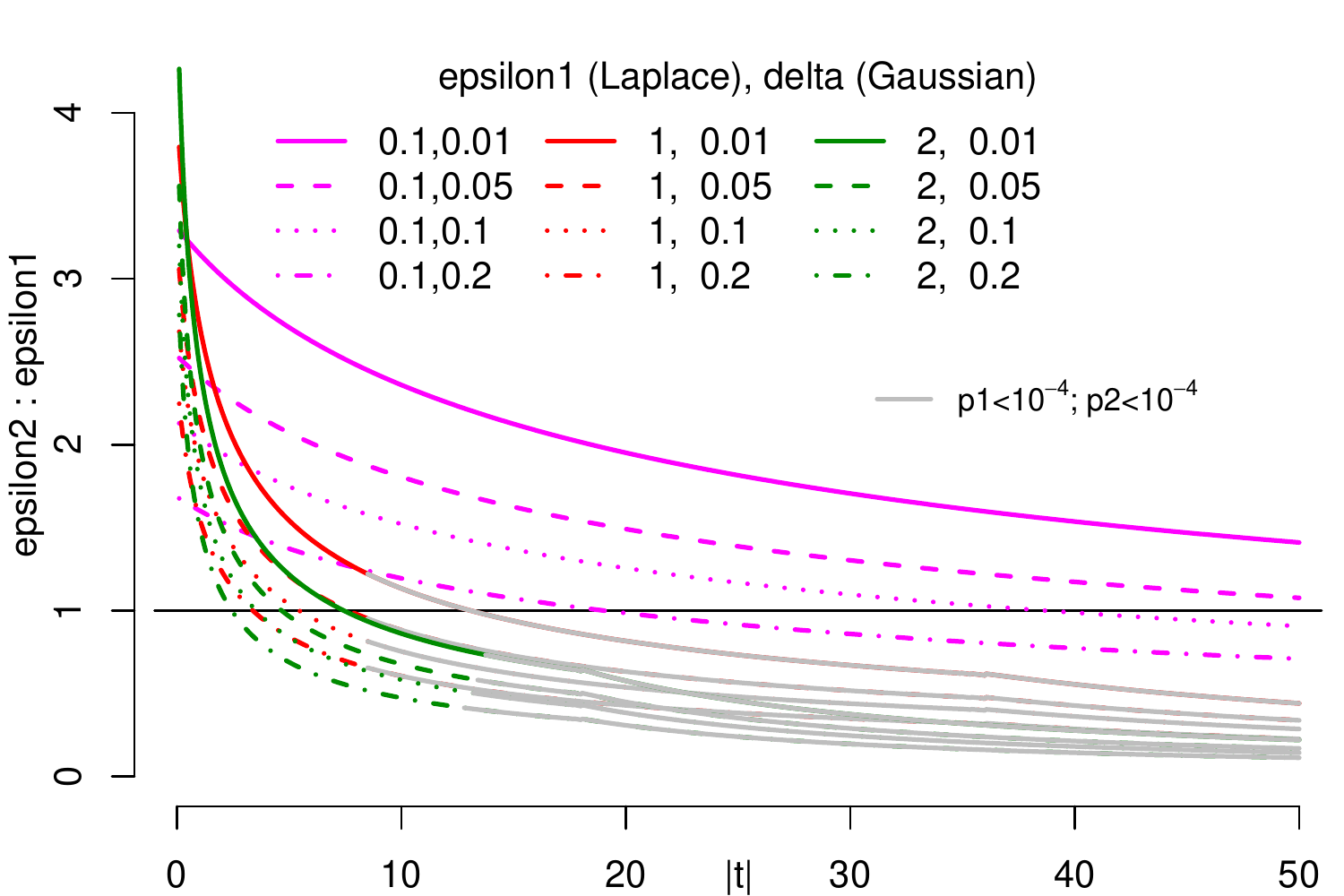}
\caption{Relative privacy cost $\epsilon_2:\epsilon_1$ (the gray curves represent the unlikely cases where both $p_1$ and $p_2$ are $<10^{-4}$)}  \label{fig:epsGaussian}
\end{center}\end{figure}

We investigate the accuracy of $\s^*$  by examining the tail probability and the dispersion of the noises injected via the $\epsilon$-DP Laplace mechanism and the $(\epsilon,\delta)$-pDP Gaussian mechanism. Denote the noise drawn from the Laplace distribution by $e_1$ and that from the Gaussian distribution by $e_2$. The location parameters of both are $\mu=0$; the tail probability $p_1=\Pr(e_1>|t|)=\exp(-|t|\epsilon/\Delta_1)$ in the Laplace distribution and $p_2=\Pr(e_2>|t|)=2\Phi(-|t|/\sigma)$ in the Gaussian distribution, where $\sigma$ is given in Eqn (\ref{eqn:lowerbound2}). Since the CDF $\Phi()$ does not have a close-formed expression, we examine several numerical examples to compare $p_1$ and $p_2$ (Figure \ref{fig:tail}). We set $\epsilon$ to be the same (0.1, 1, 2, respectively) between the two mechanisms and examine $\delta= (1\%, 5\%, 10\%,  20\%)$ for the  $(\epsilon,\delta)$-pDP Gaussian mechanism. If the ratio $p_1:p_2$ is $<1$, it implies that the Laplace mechanism is less likely to generate more extreme $\s^*$  compared to the Gaussian mechanism at the same privacy specification of $\epsilon$. We should focus on the meaningful cases where noise $|t|$ at least has a non-ignorable chance to occur in either mechanism. We used cutoff $10^{-4}$; that is, either $p_1>10^{-4}$ or $p_1>10^{-4}$ (other cutoffs can be used, depending on how ``unlikely'' is defined). It is interesting to observe that after the initial take-off at 1 when $|t|=0$, the ratio decreases until it hits the bottom and then bounds back with some cases eventually exceeding 1 at some value of $|t|$, depending on the privacy parameter specification.  The smaller $\epsilon$ or $\delta$ is, the longer it takes for the bounce-back to occurs. The observation suggests that the Laplace mechanism is in some cases more likekly to generate sanitized results $\s^*$ that are far away from $\s$.
We also compare the privacy parameter $\epsilon$ between the two mechanisms when both have the same tail probability. Figure \ref{fig:epsGaussian} shows the calculated $\epsilon_2$ value associated with the Gaussian mechanism of $(\epsilon_2,\delta)$-DP for a given $\delta$ that yields $\Pr(e_2<|t|)=\Pr(e_1<|t|)$ with the Laplace mechanism of $\epsilon_1$-DP. If the ratio of $\epsilon_2:\epsilon_1<1$  at some $|t|$ and a small and somewhat ignorable $\delta$, it implies the same tail probability can be achieved  with less privacy cost with the Gaussian mechanism compared to the Laplace mechanism. Figure \ref{fig:epsGaussian} suggests that at the same $|t|$, the more relaxation of the pure $\epsilon$-DP is allowed (i.e., the larger $\delta$ is), the smaller $\epsilon_2$ is (relative to baseline  $\epsilon_1$), which expected as the $\epsilon$ and $\delta$ together determine the noise released in the Gaussian mechanism.

Lemma \ref{lem:L1L2} presents the precision comparison of $\s^*$ between the  Laplace mechanism of $\epsilon$-DP and the Gaussian mechanism of $(\epsilon,\delta)$-pDP. With the same location parameter in the Laplace and Gaussian distributions, a larger precision is equivalent to a smaller mean squared error (MSE).
\begin{lem}\label{lem:L1L2}
Between the Gaussian mechanism of  $(\epsilon,\delta)$-pDP and the Laplace mechanism of $\epsilon$-DP for sanitizing a statistic $s$, when $\delta\!<\!2\Phi(\sqrt{2})\!\approx\!0.157$,  the variance of the Gaussian distribution in the Gaussian mechanism is always greater than the variance of the Laplace distribution associated with the  Laplace mechanism.
\end{lem}
\noindent   The proof  is provided in Appendix \ref{app:L1L2}. Lemma \ref{lem:L1L2}  suggests that  there is more dispersion in the perturbed  $s^*$ released by the Gaussian mechanism of $(\epsilon,\delta< 0.157)$-pDP than the  Laplace mechanism of $\epsilon$-DP. In other words, if there are multiple sets of $s^*$  released via the Gaussian and the Laplace mechanisms respectively, then the former sets would have a wider spread than the latter.  Since $(\epsilon,\delta)$-pDP provides less privacy protection than  $\epsilon$-pDP, together with the larger MSE, it can be argued that the Laplace mechanism is superior to the Gaussian mechanism (which is also reflected in the 3 experiments in Section \ref{sec:experiments}). It should be noted that $\delta<0.157$  in Lemma \ref{lem:L1L2} is a sufficient but not necessary condition. In other words, the Gaussian mechanism may not be less dispersed than the Laplace mechanism when  $\delta\ge0.157$.  Furthermore,  since $\delta$ needs to be small to provide sufficient privacy protection in the setting of $(\epsilon,\delta)$-pDP, it is  very unlikely to have $\delta>0.157$ in practical applications.    Also noted is that the setting explored in Lemma \ref{lem:L1L2}, where the focus is on examining the precision (dispersion) of a single perturbed  statistic given the specificized privacy parameters and the original statistics when the sample size of a data set is public, is different from the recent work on the bounds of sample complexity (required sample size) to reach a certain level of a statistical \emph{accuracy}  in perturbed  results with $\epsilon$-DP or  $(\epsilon,\delta)$-aDP \cite{harvard} (more discussions are provided in Section \ref{sec:discussion} on this point).

\section{Experiments}\label{sec:experiments}
We run three experiments on the mildew data set, the Czech data set, and the Census Income data set; a.k.a. the adult data.    The mildew data contains information of parental alleles at 6 loci on the chromosome for 70 strands of barley powder mildew\cite{charest}. Each loci has two levels, yielding  a very sparse 6-way cross-tabulation (22 cells out of the 64 are non-empty with low frequencies in many other cells). The Czech data contains data collected on 6 potential risk factors for coronary thrombosis for 1841 workers in a Czechoslovakian car factory \cite{charest}. Each risk factor has 2 levels (Y or N).  The cross-tabulation is also 6-way  with 64 cells, the same as the mildew data, but table is not as sparse with the large $n$ (only one empty cell). The adult data was extracted from the 1994 US Census database to yield a set of reasonably clean records that satisfy a set of conditions\cite{adult}. The data set is often used to test classifiers by predicting  whether a person makes over 50K a year. We used only the completers in the adult data (with no missing values on the attributes) and then split them to 2/3 training (20009 subjects) and 1/3 testing (10005 subjects).
\begin{figure}[!htb] \begin{center}
\includegraphics[scale=0.55]{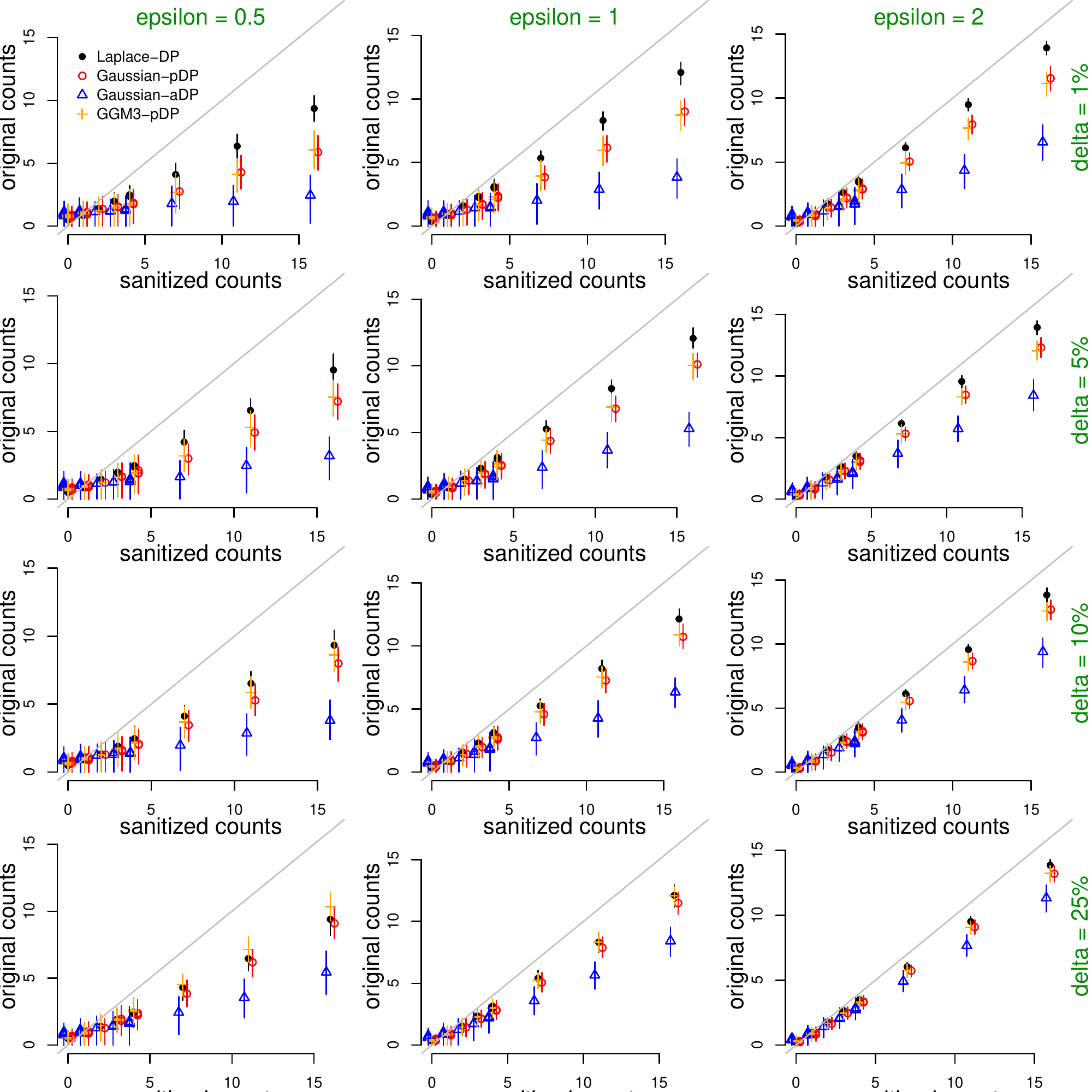}
\caption{sanitized vs. original cell counts in the mildew data}\label{fig:mildew.count}
\end{center}\end{figure}
\begin{figure}[!htb] \begin{center}
\includegraphics[scale=0.6]{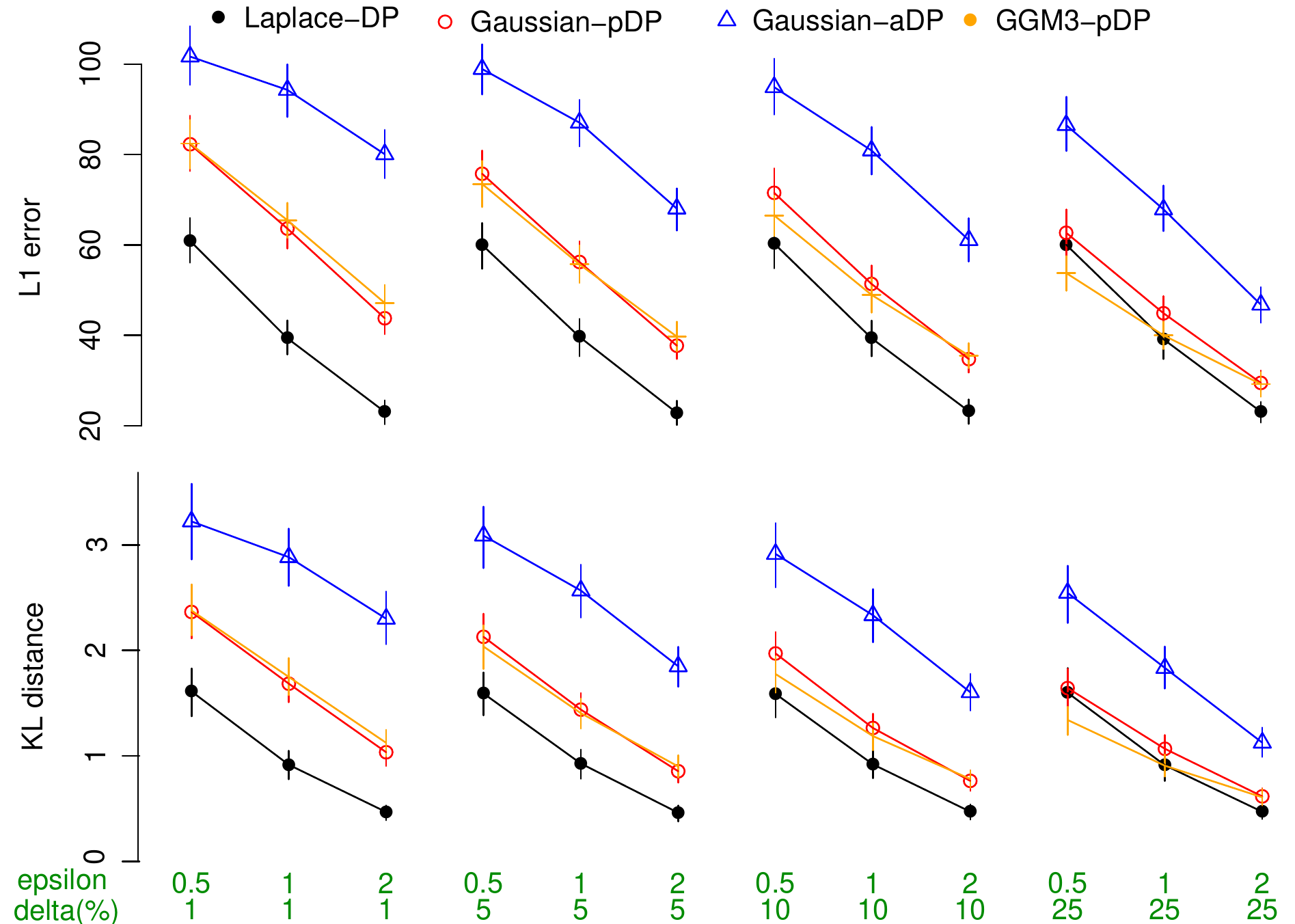}
\caption{$l_1$ distance and KL divergence between sanitized and original counts in the mildew data}\label{fig:mildew.L1KL}
\end{center}\end{figure}
\begin{figure}[!htb] \begin{center}
\includegraphics[scale=0.55]{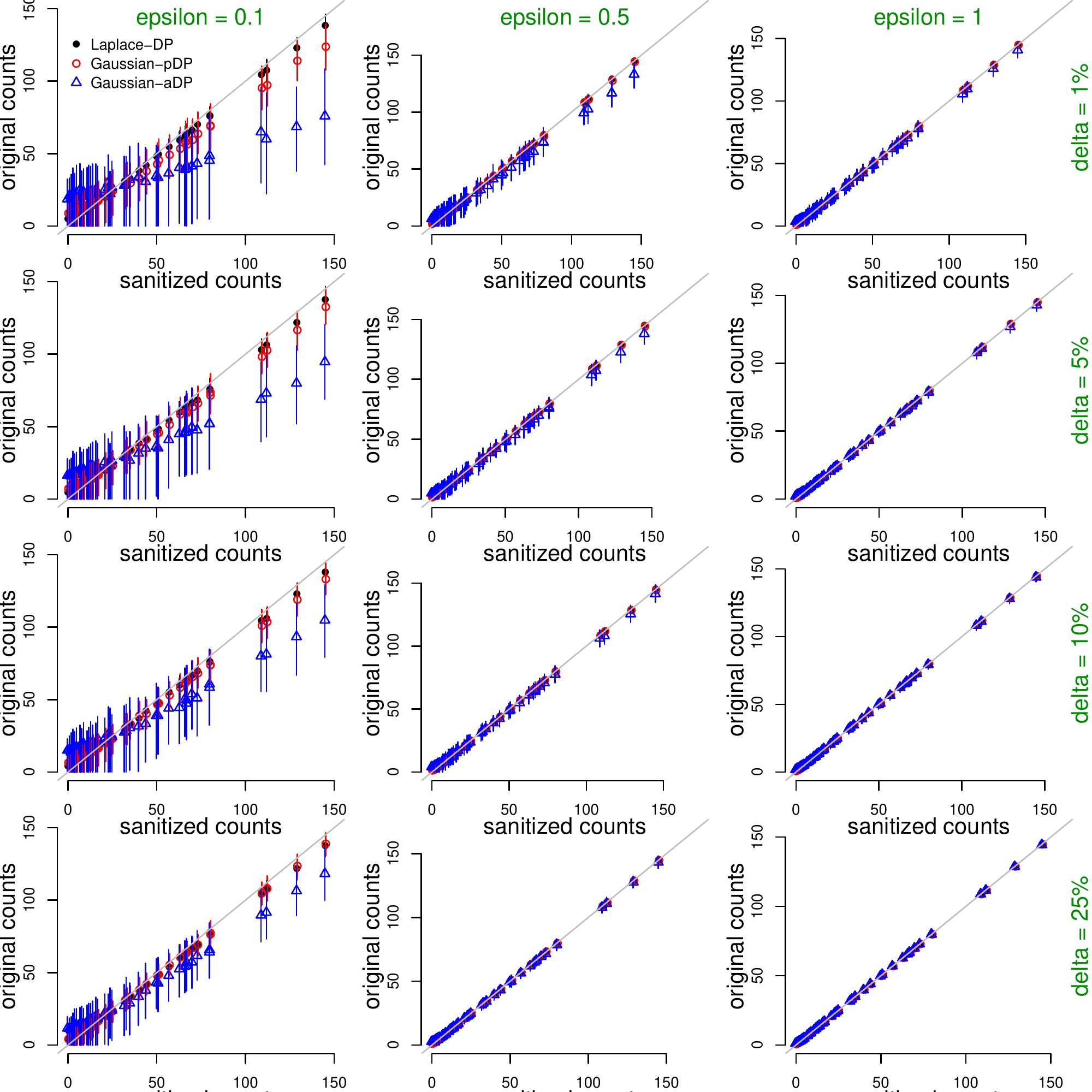}
\caption{sanitized vs. original cell counts in the Czech data}\label{fig:czech.count}
\end{center}\end{figure}
\begin{figure}[!htb] \begin{center}
\includegraphics[scale=0.6]{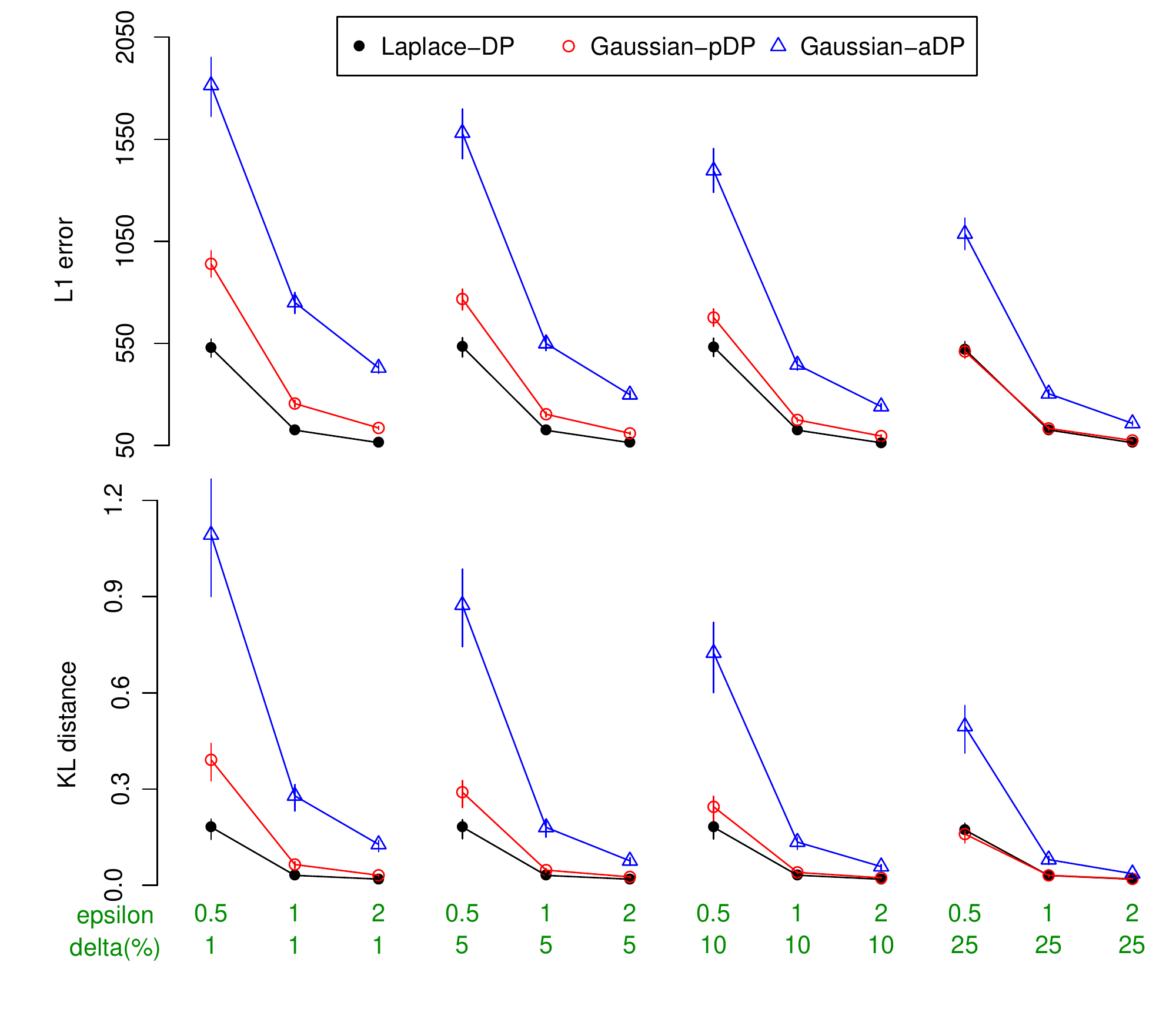}
\caption{$l_1$ distance and KL divergence between sanitized and original counts  in the Czech data}\label{fig:czech.L1KL}
\end{center}\end{figure}

In each experiment, we run the Laplace mechanism of $\epsilon$-DP, the Gaussian mechanism of $(\epsilon,\delta)$-pDP presented in Section \ref{sec:gaussian}, and the Gaussian mechanism of of $(\epsilon,\delta)$-aDP \cite{privacybook} to sanitize count data. We examined $\epsilon=0.5, 1,2$ and $\delta=0.01, 0.05, 0.1, 0.25$. To examine the variation of noises, we run 500 repeats and computed the means and standard deviations of $l_1$ distances  between the sanitized and the original counts and the Kullback-Leibler (KL) divergence between the empirical distributions of the synthetic data and the original data over the 500 repeats.  In addition, we  tested the GG mechanism of order 3 ($p=3$) in the mildew data, and compared the classification accuracy of the income outcome in the testing data set in the adult experiment based on the support vector machines (SVMs) trained with the original training data and the sanitized training data, respectively.   The KL distance was calculated using the \texttt{KL.Dirichlet} command in R package \texttt{entropy} that computes a Bayesian estimate of the KL divergence. The SVMs were trained using the \texttt{svm} command in R package \texttt{e1071}. In all experiments,  $\Delta_p=1$ for all $p$  since the released query is a histogram and the bin counts are based on disjoint subsets of data. The scale parameters of the Laplace mechanism and the Gaussian mechanisms were obtained analytically  ($\Delta_1\epsilon^{-1}$, Eqns (\ref{eqn:lowerbound2}) and (\ref{eqn:privacybook}), respectively), the grid search and the MC approach were applied to obtain the lower bound $b$ for GGM-3 via Corollary \ref{cor:disjoint}. In the mildew and Czech experiments, we sanitized all bins in the histograms, including the empty bins, assuming all combinations of the 6 attributes in each case are practically meaningful (in other words, the empty cells are sample zeros rather than population zeros). In the adult data, there are 14  attributes and $\sim\!1.944\times10^{13}$ bins in the 14-attribute histogram, a non-ignorable portion of which do not make any practical sense (e.g., a 90-age works $>80$ hours per week). For simplicity, we only sanitized the 17,985 nonempty cells in the training data.  After the sanitization, we set the out-of-bounds synthetic counts $<0$ at 0 and those $>n$ at  $n$, respectively, and normalized the sanitized counts to sum up to the original sample size $n$ in all 3 experiments, assuming $n$ itself is public or does not carry privacy information.
\begin{figure}[!htb] \begin{center}
\includegraphics[scale=0.55]{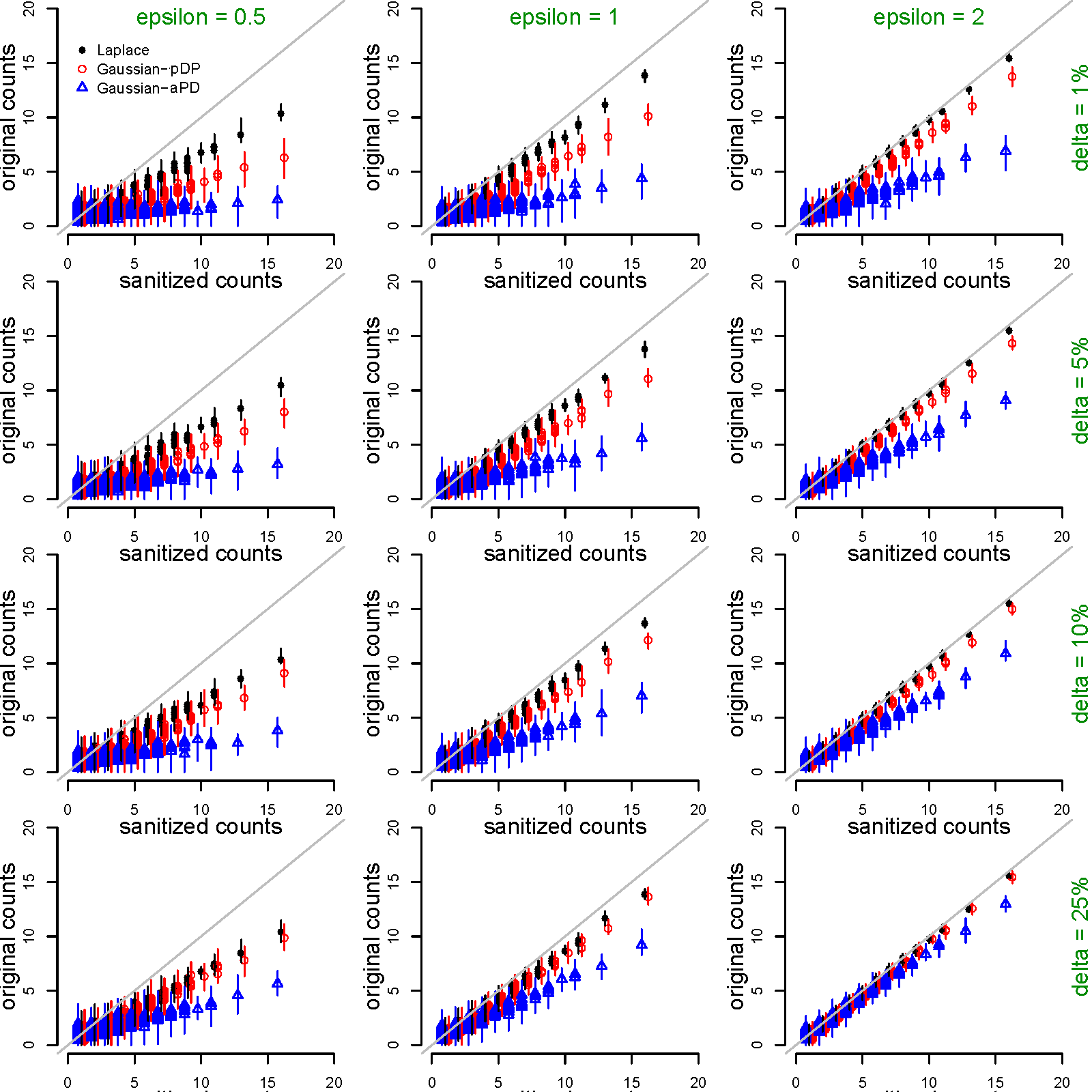}
\caption{sanitized vs. original cell counts in the adult data}\label{fig:adult.count}
\end{center}\end{figure}
\begin{figure}[!htb] \begin{center}
\includegraphics[scale=0.6]{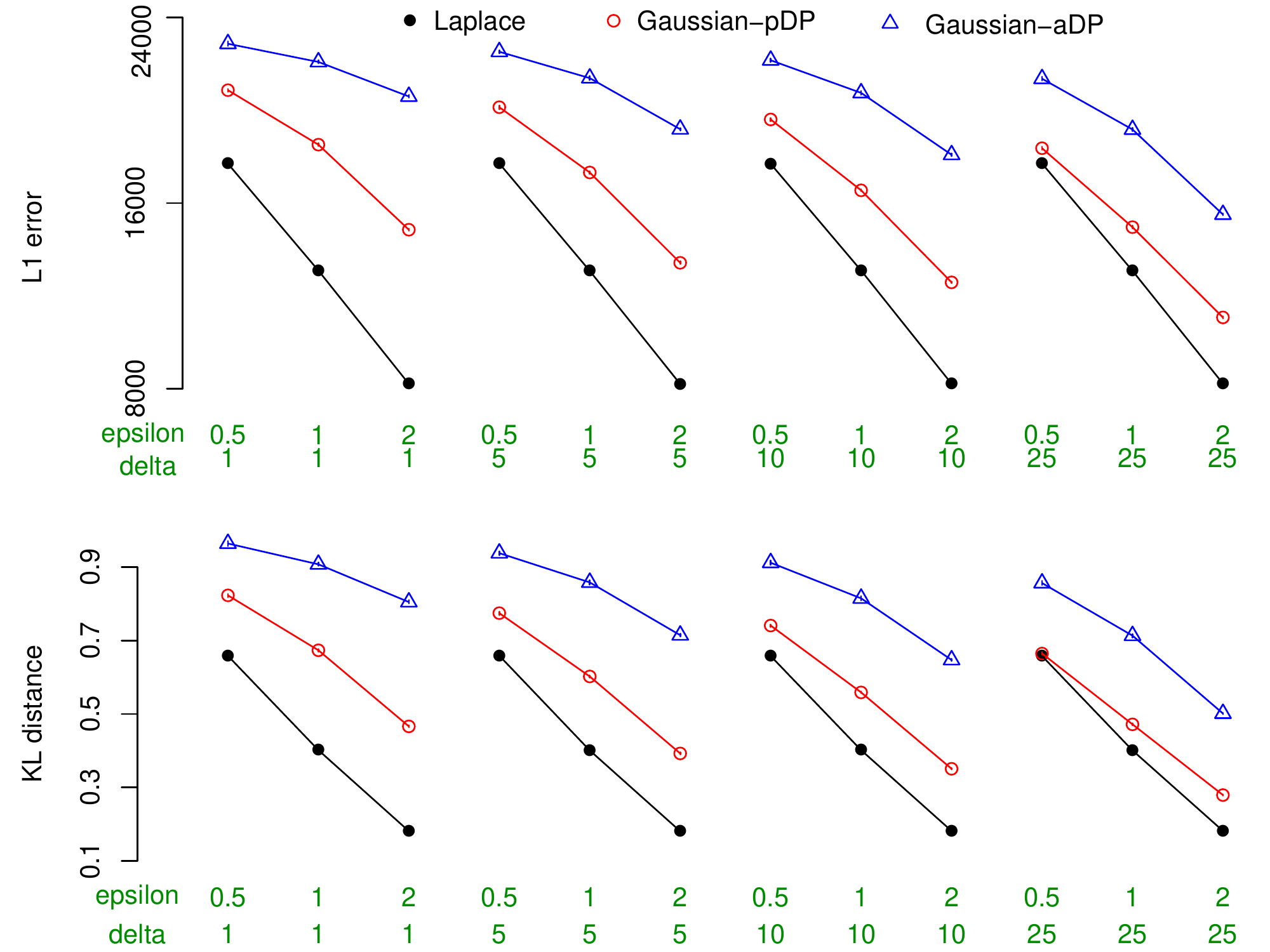}
\caption{$l_1$ distance and KL divergence between sanitized and original counts in the adult data}  \label{fig:adult.L1KL}
\end{center}\end{figure}
\begin{figure}[!htb] \begin{center}
\includegraphics[scale=0.6]{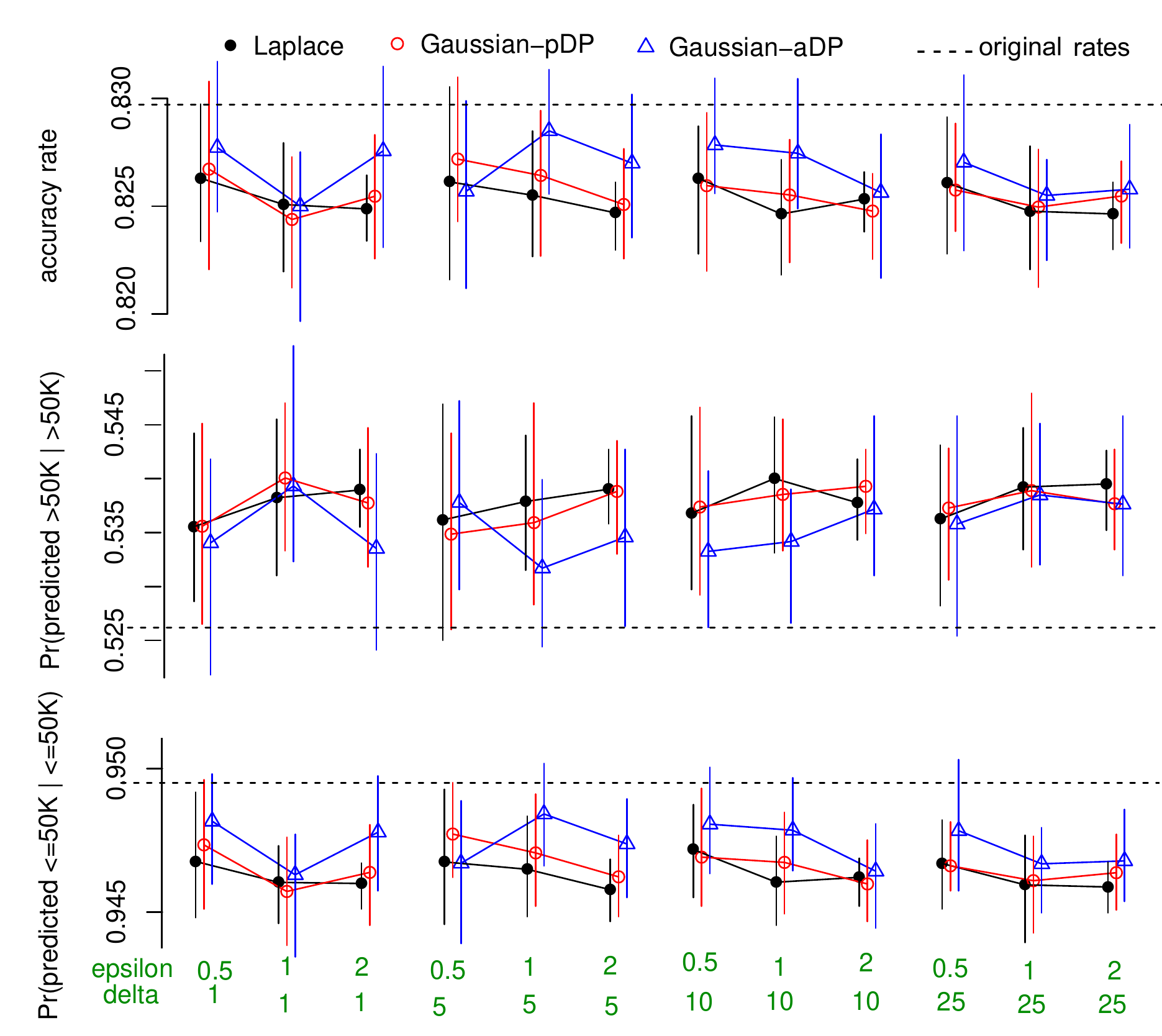}
\caption{Prediction accuracy in testing data via SVMs trained on sanitized and original data in the adult data}\label{fig:adult.acc}
\end{center}\end{figure}

The results are given in Figures \ref{fig:mildew.count} to \ref{fig:adult.acc}. In Figures \ref{fig:mildew.count}, \ref{fig:czech.count} and \ref{fig:adult.count}, the closer the points are to the identity line, the more similar are the original and sanitized counts.  The Laplace sanitizer is the obvious winner in all 3 cases, producing the sanitized counts closest to the original with the smallest $l_l$ error and the KL divergence, followed by the Gaussian mechanism of  $(\epsilon,\delta)$-pDP, and GGM3  of  $(\epsilon,\delta)$-pDP in the mildew data; the Gaussian mechanism of $(\epsilon,\delta)$-aDP is the worst. In the mildew experiment, the performance of the Gaussian mechanism of $(\epsilon,\delta)$-pDP is similar when $\epsilon=2$ or $\delta\ge0.1$. The decrease in the $l_1$ error and the KL divergence seems to decrease more or less in a linear manner as $\epsilon$ increases from 0.5 to 1 to 2, while the impact of $\delta$ seemed to have less  a profound impact on the  $l_1$ error and the KL divergence. In the Czech experiment,  the sanitized counts approach the original counts more quickly than the mildew case with increased $\epsilon$ and $\delta$, but there is significantly more variability for small $\epsilon$ (0.1); and the  $l_1$ error and the KL divergence no longer decreases in a linear fashion, but drastically from $\epsilon=0.5$ to 1 and much less from $\epsilon=1$ to 2. The differences in the  results between the mildew and the Czech experiments can be explained by the larger $n$ in the latter. In the adult experiment,  Figure \ref{fig:adult.acc} suggests the prediction accuracy via the SVMs built on sanitized data is barely affected compared to the original accuracy regardless of the mechanism.There are some decreases in the accuracy rates from the original, but they are largely ignorable (on the scale of 0.25\% to 1\%), even with the variation take into account. In addition, the Gaussian mechanism of $(\epsilon,\delta)$-aDP, though being the worst in preserving the original counts measured the $l_1$ distance and KL divergence, is no worse than the two Gaussian mechanisms in prediction.

\section{Discussion}\label{sec:discussion}
We introduced a new concept of the $l_p$ GS, and unified the Laplace mechanism and the Gaussian mechanism in the family of the GG mechanism. For bounded data, we discussed the truncated and the BIT GG mechanisms to achieve $\epsilon$-DP.  We also proposed $(\epsilon,\delta)$-pDP as an alternative paradigm to the pure $\epsilon$-DP for the GG mechanism  for order $p\ge2$. We showed the connections and distinctions between the GG mechanism and the Exponential mechanism when the utility function is defined as the negative $p\tthh$-power  of the Minkowski distance between the original and  sanitized results. We also presented the Gaussian mechanism as an example of the GG mechanism and derived a lower bound for the scale parameter of the associated Gaussian distribution to achieve $(\epsilon,\delta)$-pDP. The bound is tighter than the lower bound for the Gaussian mechanism of $(\epsilon,\delta)$-aDP. We compared the tail probability and the dispersion of the the noise generated via the  Gaussian mechanism of $(\epsilon,\delta)$-pDP and the Laplace mechanism. We finally applied the Gaussian mechanisms of $(\epsilon,\delta)$-pDP and $(\epsilon,\delta)$-aDP and the Laplace mechanism of $\epsilon$-DP in three real-life data sets.

The GG mechanism is based on the  $l_p$ $``$global$"$ sensitivity of query results in the sense that the sensitivity is independent of any specific data. Though the employment of the GS  is robust  in terms of privacy protection, it could result in a large amount of noises being injected to query results. There is work that allows the sensitivity of a query to vary with data  ($``$local$"$  sensitivity) \cite{smooth, robust} with the purpose to increase the accuracy of sanitized  results. How to develop the GG mechanism in the context of local sensitivity is a topic for future investigation.

The setting for the examination on the tail probability and dispersion in Section \ref{sec:gaussian} is  different from, though  related to,  the work on upper and lower bounds on sample complexity -- the required sample size $n$ to reach a certain level of accuracy  $\alpha$ and  privacy guarantee $(\epsilon,\delta)$ for  count queries \cite{geometry, fingerprinting, harvard}.   $\alpha$ often refers to the accuracy of perturbed  results in the DP literature, such as  the worst case accuracy $L_\infty$ or average accuracy  $L_1,$ and might also refer to the tail probability and the MSE of released data, among others. A differential privacy mechanism is characterized by  $\epsilon$ (and $\delta$) for privacy guarantee, $\alpha$ to measure information preservation and utility of sanitized  results, and the sample size $n$ of original data.
The existing work on sample complexity focuses on bounding $n$ given $\epsilon$ (and $\delta)$ and $\alpha$, while the results in  Section \ref{sec:gaussian} focus on the the accuracy and precision of sanitized results given $\epsilon$ (and $\delta)$ and $n$. If the bias from perturbed  results (relative to the original results) are the same between the two mechanisms, a larger precision is equivalent to a smaller MSE.

\section*{Appendix}
\begin{appendix}
\renewcommand{\theequation}{\Alph{section}.\arabic{equation}}
\setcounter{equation}{0}
\section{\normalsize{Proof of Lemma \ref{lem:GSp1}}}\label{app:GSp1}
\vspace{-12pt}
\begin{align*}
\Delta_{p}&\!=\!\textstyle\maxxtex\!\!  \left(\sum_{k=1}^r\left|s_k(\x)\!-\!s_k(\x')\right|^p\right)^{\!1/p}
=\!\left(\textstyle\maxxtex\! \! \sum_{k=1}^r\left|s_k(\x)\!-\!s_k(\x')\right|^p\!\right)^{\!1/p}\!\!\!\!.\\
&\mbox{\hspace{12pt}Since }\maxxtex \textstyle\sum_{k=1}^r\left|s_k(\x)\!-\!s_k(\x')\right|^p
\le\textstyle\sum_{k=1}^r\!\maxxtex\! \! \left|s_k(\x)-s_k(\x')\right|^p\\
&=\textstyle\sum_{k=1}^r\! \! \left(\maxxtex \left|s_k(\x)-s_k(\x')\right|\right)^p
=\textstyle\sum_{k=1}^r\Delta_{1,k}^p.
\end{align*}
Therefore, $\left(\!\sum_{k=1}^r\!\Delta_{1,k}^p\!\right)^{\!1/p}$ is an upper bound for $\Delta_{p}.\blacksquare$

\section{\normalsize{Proof of Claim \ref{cla:lowerb}}}\label{app:GGM}
\allowdisplaybreaks
\begin{align}
&\left|\log\!\!\left(\!\frac{\Pr(\s^{\ast} \in Q|\x)}{\Pr(\s^{\ast}\in Q|\x')}\!\right)\right|
=\left|\log\!\!\left(\frac{ \exp\left(-b^{-p}\|\s^*-\s(\x)\|_p^p\right) }{\exp\left(-b^{-p}\|\s^*-\s(\x')\|_p^p\right) }\!\right)\right|\notag\\
&=b^{-p}\big| \|\s^*-\s(\x)\|_p^p- \|\s^*-\s(\x')\|_p^p\big|=b^{-p}\big|\textstyle\sum_{k=1}^r\left(|s_k^*-s_k(\x)|^p- |s_k^*-s_k(\x')|^p\right)\big|\notag\\
&\le b^{-p}\textstyle\sum_{k=1}^r\big||s_k^*-s_k(\x)|^p- |s_k^*-s_k(\x')|^p\big|
=b^{-p}\textstyle\sum_{k=1}^r\!\big||e_k|^p- |e_k+d_k|^p\big|,\label{eqn:first} \\
& \mbox{ where  $e_k=s_k^*-s_k(\x)$ and $d_k=s_k(\x)-s_k(\x')$}\notag \\
&=b^{-p}\textstyle\sum_{k=1}^r\!\big||e_k^p|- |(e_k+d_k)^p|\big|\mbox{ for integers $p\ge1$}\notag\\
&\le b^{-p}\textstyle\sum_{k=1}^r\!\big|e_k^p-(e_k+d_k)^p\big|=b^{-p}\textstyle\sum_{k=1}^r\!\left|\sum_{j=1}^p (_j^p)e_k^{p-j}d_k^j\right| \mbox{ by reverse triangle inequality} \notag\\ 
& \le b^{-p}\textstyle\sum_{k=1}^r\sum_{j=1}^p(_j^p)|e_k|^{p-j}|d_k|^j\label{eqn:second}\\
&= b^{-p}\!\left(\!\textstyle p\sum_{k=1}^r\!|e_k|^{p-1}|d_k|\!+\!\frac{(p-1)p}{2}\sum_{k=1}^r\!|e_k|^{p-2} |d_k|^2\!+\cdots\!+\! (p-1)p\textstyle\sum_{k=1}^r\!|e_k|^2|d_k|^{p-2}/2\!\right.\notag\\
&\qquad+\textstyle\left.\! p\! \sum_{k=1}^r\!|e_k|\cdot|d_k|^{p-1} + \!\sum_{k=1}^r|d_k|^p \right)\notag\\
\!\!&\!\!\!\le\! b^{-p}\textstyle\left(p\!\sum_{k=1}^r\!|e_k|^{p-1}\!\Delta_{1,k}\! +\!\frac{(p-1)p}{2}\sum_{k=1}^r\!|e_k|^{p-2}\! \Delta_{1,k}^2\!+\cdots\right.\notag\\
&\qquad\quad\textstyle\left. \!+ \frac{(p-1)p}{2}\sum_{k=1}^r\!|e_k|^2 \Delta_{1,k}^{p-2}\!+\!p\! \sum_{k=1}^r\!|e_k|\Delta_{1,k}^{p-1}\!+\!\Delta_{p}^p\right),\label{eqn:third}
\end{align}
where $\Delta_{1,k}$ is the $l_1$ GS of $s_k$ and $\Delta_{p}$ is the $l_p$ GS of $\s$. To achieve $\epsilon$-DP, Eqn (\ref{eqn:third}) needs to be $\le\epsilon$; that is,
\begin{equation}\label{eqn:fourth}
\textstyle\Delta_{p}^p+\sum_{k=1}^r\sum_{j=1}^{p-1}(_j^p)|e_k|^{p-j}\Delta_{1,k}^j\le b^p\epsilon.
\end{equation}
A less tight bound can be obtained by applying Lemma \ref{lem:GSp1} ($\Delta_{p}^p\le\sum_{k=1}^r\Delta_{1,k}^p$), thus
\begin{equation}\label{eqn:fifth}
\textstyle\sum_{k=1}^r\sum_{j=1}^p(_j^p)|e_k|^{p-j}\Delta_{1,k}^j\le b^p\epsilon.
\end{equation}
The inequalities in Eqns (\ref{eqn:fourth}) or (\ref{eqn:fifth}) susgest that the lower bound on $b$ depends on the random GG noise $e_k=s^*_k-s_k$ for $k=1,\ldots,r$, the support of which is $(-\infty,\infty)^r$. In other words, there does not exist a random noise-free solution on $b$, unless $p=1$ in which case the inequality no longer involves the error terms and the GG mechanism reduces to the familiar Laplace mechanism of $\epsilon$-DP, leading to Claim \ref{cla:lowerb}. When $p=1$,  Eqn (\ref{eqn:first}) $\!\!\le b^{-1}\!\sum_{k=1}^r\!|d_k|\!\le\! b^{-1}\textstyle\sum_{k=1}^r\!|\Delta_{1,k}|\!=\!b^{-1}\Delta_{1}\!<\!\epsilon$, and thus $b\!>\!\Delta_{1}\epsilon^{-1}$. $\blacksquare$

\section{\normalsize{Proof of  \texorpdfstring{$\epsilon$}{}-DP of the truncated GG mechanism in Definition \ref{def:tGGM}}}\label{app:tGGM}
To satisfy $\epsilon$-DP, we need
\begin{align}
&\left|\log\!\left(\!\frac{\Pr(\s^{\ast} \in Q|\x, \s^{\ast}\in[c_{10},c_{11}]\times\!\cdots\!\times[c_{r0},c_{r1}])}{\Pr(\s^{\ast}\in Q|\x', \s^{\ast}\in[c_{10},c_{11}]\times\!\cdots\!\times[c_{r0},c_{r1}])}\!\right)\right|\notag\\
=&\left|\log\!\left(\frac{ \exp\left(-b^{-p}\|\s^*-\s(\x)\|_p^p\right)}{\prod_{k=1}^r\Pr(c_{k0}\! \le \!s^*_k\! \le\! c_{k1};  s_k, b,p)}\times\frac{\prod_{k=1}^r\Pr(c_{k0}\! \le \!s^*_k\! \le\! c_{k1};  s'_k, b,p)}{\exp\left(-b^{-p}\|\s^*-\s(\x')\|_p^p\right) }\!\right)\right|\notag\\
=&\left|\log\!\left(\!\frac{ \exp\left(-b^{-p}\|\s^*-\s(\x)\|_p^p\right) }{\exp\left(-b^{-p}\|\s^*-\s(\x')\|_p^p\right) }\!\right)+\log\!\left(\frac{ \prod_{k=1}^r\Pr(c_{k0}\! \le\!s^*_k\! \le\! c_{k1};  s_k, b,p)}{\prod_{k=1}^r\Pr(c_{k0}\! \le\!s^*_k\!\le\! c_{k1};  s'_k, b,p)}\right)\right|\notag\\
\le&\left|\log\!\left(\!\frac{ \exp\left(-b^{-p}\|\s^*-\s(\x)\|_p^p\right) }{\exp\left(-b^{-p}\|\s^*-\s(\x')\|_p^p\right) }\!\right)\right|+\label{eqn:tGGD1}\\
&\left|\log\!\left(\frac{ \prod_{k=1}^r\Pr(c_{k0}\! \le\!s^*_k\! \le\! c_{k1};  s_k, b,p)}{\prod_{k=1}^r\Pr(c_{k0}\! \le\!s^*_k\!\le\! c_{k1};  s'_k, b,p)}\right)\right|\le\epsilon\label{eqn:tGGD2}
\end{align}
If  the term in Eqn (\ref{eqn:tGGD1}) satisfies $\epsilon/2$-DP, so does Eqn (\ref{eqn:tGGD2}). Appendix \ref{app:GGM} establishes that Eqn (\ref{eqn:tGGD2}) satisfies $\epsilon/2$-DP when $b^p(\epsilon/2)\ge \Delta_{p}^p+\sum_{k=1}^r\sum_{j=1}^{p-1}(_j^p)|s^*_k-s_k|^{p-j}\Delta_{1,k}^j$ Since $\s^*$ is bounded within $[c_{k0},c_{k1}]$  for $k=1,\ldots,K$, $|s^*_k-s_k|\le|c_{k1}-c_{k0}|$. Setting $b^p(\epsilon/2)\ge \Delta_{p}^p+\sum_{k=1}^r\sum_{j=1}^{p-1}(_j^p)|c_{k1}-c_{k0}|^{p-j}\Delta_{1,k}^j$ ensures the truncated GG mechanism is of $\epsilon$-DP; or equivalently, $b^p \ge 2\epsilon^{-1}\!\!\textstyle\left(\!\sum_{k=1}^r\sum_{j=1}^{p-1}(_j^p)|c_{k1}-c_{k0}|^{p-j}\Delta_{s_k}^j\!+\!\Delta_{\s,p}^p\!\right)$  ensures that the truncated GG mechanism is of $\epsilon$-DP. $\quad\blacksquare$

\section{\normalsize{Conservativeness of Exponential mechanism}}\label{app:exp}
\begin{cor}\label{lem:conservative}
The actual privacy cost of the Exponential mechanism of $\epsilon$-DP is always less than the nominal budget $\epsilon$. When the normalization factor $A(\x)$ in Eqn (\ref{eqn:exp}) is independent of $\x$, the actual privacy cost is $\epsilon/2$.
\end{cor}
$A(\x)$ independent of $\x$ implies increases and decreases in the utility scores upon the change  from $\x$ to $\x'$ $``$cancel out$"$ when integrated or summed over all possible $\s^*$ in the form of $\exp\!\left(u(\s^{\ast}|\x)\frac{\epsilon}{2\Delta_u}\right)$.
\begin{proof}
Since $ u(\s^{\ast}|\x)-u(\s^{\ast}|\x')\le\Delta_u$,
\begin{align}
&\left|\log\!\left(\!\frac{\Pr(\s^{\ast}(\x) \in Q)}{\Pr(\s^{\ast}(\x')\in Q)}\!\right)\right|\!=\!\left|\log\!\!\left(\!\frac{\exp\left(u(\s^{\ast}|\x)
\frac{\epsilon}{2\Delta_u}\!\right)}{\exp\left(u(\s^{\ast}|\x')\frac{\epsilon}{2\Delta_u}\right)}\!\!\times\!\frac{A(\x')}{A(\x)}\!\right)\!\right|
\le \!\left|\log\!\left(\!e^{\epsilon/2}\frac{A(\x') }{A(\x)}\!\right)\!\right| \label{eqn:ratio1}\\
&=\left|\frac{\epsilon}{2}+ \log\left(\frac{A(\x') }{A(\x)} \right)\right|
\le \frac{\epsilon}{2}+ \left|\log\left(\frac{A(\x')}{A(\x)}\right)\right|\label{eqn:ratio2}
\end{align}
by the triangle inequality, and
\begin{align}
A(\x') =& \int_{\s^*\in\mathcal{S}}\exp\!\left(u(\s^{\ast}|\x')\frac{\epsilon}{2\Delta_u}\right) d\s^{\ast}  \le\int_{\s^*\in\mathcal{S}}\exp\!\left((u(\s^{\ast}|\x)+\Delta_u)\frac{\epsilon}{2\Delta_u}\right) d\s^{\ast} \label{eqn:ratio3}\\
=&\exp\!\left(\frac{\epsilon}{2}\right)\!\!\int_{\s^*\in\mathcal{S}}\exp\!\left(u(\s^{\ast}|\x)\right)d\s^{\ast} = \exp\!\left(\frac{\epsilon}{2}\right)A(\x)\notag
\end{align}
Therefore, $\log\left(\frac{A(\x')}{A(\x)}\right)\le\epsilon/2$,  and Eqn (\ref{eqn:ratio2}) becomes
\begin{align}
\!\!\left|\log\!\!\left(\frac{\Pr(\s^{\ast}(\x) \in Q)}{\Pr(\s^{\ast}(\x')\in Q)}\right)\right|
\le \frac{\epsilon}{2}\!+\!\left|\log\left(\frac{A(\x)}{A(\x')} \right)\right|\!\le\!\epsilon\label{eqn:ratio4}
\end{align}
\normalsize  The same result can be obtained by replacing the  integral with summation when $\mathcal{S}$ is a discrete set in the equation set (\ref{eqn:ratio4}). The above results seem to suggest $\epsilon$ can be achieved exactly since $``$equality$"$ appears in all the inequalities above (Eqn  (\ref{eqn:ratio1}) to  (\ref{eqn:ratio4})); however, equality  cannot occur simultaneously in Eqns (\ref{eqn:ratio1}) and (\ref{eqn:ratio3}) unless $\Delta_u$ was 0, which is meaningless in DP.  In addition, $\Delta_u$ is defined as the maximum change in $u$ for all $d(\x,\x')=1$. While it is likely that the maximum change occurs at more than a single value of $\s^*$, it is not possible that the utility scores at all values of $\s^*$ increase or decreases by the same amount $\Delta_{u}$. In other words, the $``$equality$"$ in Eqn (\ref{eqn:ratio3}) itself is unlikely to hold. All taken together, the actual privacy cost in the Exponential mechanism is always less than $\epsilon$ and never attains the exact upper bound $\epsilon$. In the extreme, the actual privacy cost can be down to $\epsilon/2$  when $A(\x)\equiv A(\x')\; \forall\; \x,\x'$ and $d(\x,\x')=1$, as suggested by Eqn (\ref{eqn:ratio2}).
\end{proof}

\section{\normalsize{Proof of Lemma \ref{lem:deltau.s.relationship}}}\label{app:deltau.s.relationship}
\begin{proof}
\noindent \textbf{Part a)}. Denote $\s(\x)$ by $\s$ and $\s(\x')$ by $\s'$. When $p=1$, $u(\s^*|\x)=-\|\s^{\ast}-\s\|_1$, $|u(\s^{\ast}|\x)-u(\s^{\ast}|\x')|
=\!\big|\! \sum_{k=1}^r\!(|s^{\ast}_k-s_k|-|s^{\ast}_k-s'_k|)\big|
\!\le\!  \sum_{k=1}^r \!\big| |s^{\ast}_k-s_k|-|s^{\ast}_k-s'_k|\big|
\!\le\! \sum_{k=1}^r \!\big| s^{\ast}_k-s_k- (s^{\ast}_k-s_k)\big|
\!= \! \sum_{k=1}^r\! |s_k-s'_k|\!=\!|\s-\s'|_1$. Therefore, $\Delta_u=\!\maxxutex\!|u(\s^{\ast}|\x)-u(\s^{\ast}|\x')|\!\le\!\maxxtex\!\|\s-\s'\|_1=\Delta_{\s,1}$.
\end{proof}
\begin{proof}
\noindent \textbf{Part b)}. When $p=2$, $u(\s^*|\x)=-\|\s^{\ast}-\s\|^2_2$, $|u(\s^{\ast}|\x)-u(\s^{\ast}|\x')|=\big|\sum_{k=1}^r(s_k-s_k^{\ast})^2-(s'_k-s_k^{\ast})^2\big|\le
\sum_{k=1}^r \big|(s_k-s_k^{\ast})^2-(s'_k-s_k^{\ast})^2\big|=
\sum_{k=1}^r |s_k-s'_k|\cdot|s_k-s_k^{\ast}+s'_k-s_k^{\ast}|\le
\sum_{k=1}^r\Delta_{1,k}(|s_k-s_k^{\ast}|+|s'_k-s_k^{\ast}|)$. Suppose $s_k$ is bounded within $[c_{k0},c_{k1}]$, so is $s_k^{\ast}$, then
\begin{align}\label{eqn:deltau2}
\Delta_u&=\!\!\!\!\textstyle\maxxu \big|\sum_{k=1}^r (s_k(\x)-s_k^{\ast})^2\!-\!\sum_{k=1}^r (s_k(\x')-s_k^{\ast})^2\big|\le2\textstyle\sum_{k=1}^r\Delta_{1,k}(c_{k1}-c_{k0})
\end{align}
When $c_{k1}-c_{k0}\equiv b-a\;\forall\;k$, $\Delta_u\le2(b-a)\sum_{k=1}^r\Delta_{1,k} =2(b-a)\Delta_{1}$.
\end{proof}
\begin{proof}
\noindent \textbf{Part c)}. When $u(\s^*|\x)=-\|\s^{\ast}-\s\|^p_p$ for integer $p\ge1$,  $|u(\s^{\ast}|\x)-u(\s^{\ast}|\x')|= \big|\|\s^{\ast}-\s\|^p_p-\|\s^{\ast}-\s'\|^p_p\big|=\big|\sum_{k=1}^r|s_k-s_k^{\ast}|^p-\sum_{k=1}^r|s'_k-s_k^{\ast}|^p\big|\le\sum_{k=1}^r\big\|(s_k-s_k^{\ast})^p|-|(s'_k-s_k^{\ast})^p|\big|
\le\sum_{k=1}^r\big|(s_k-s_k^{\ast})^p-(s'_k-s_k^{\ast})^p\big|
\!=\!\sum_{k=1}^r\big|\sum_{i=1}^p(^p_i)(-s_k^{\ast})^{p-i}\left[ s_k^i -(s'_k)^i\right]\big|
\le\sum_{k=1}^r\sum_{i=1}^p (^p_i)\big|(s_k^{\ast})^{p-i}\left[ s_k^i -(s'_k)^i\right]\big|$. Suppose $s_k$ is bounded within $(c_{k0},c_{k1})$, so is $s_k^{\ast}$. \\ Define $\Delta_{1,k}^{(i)}=\maxxtex|s_k^i -(s'_k)^i|$ , then
\begin{align}\label{eqn:deltaup}
&\Delta_u\!=\!\maxxutex\big|\!\textstyle\sum_{k=1}^r\! |s_k-s_k^{\ast}|^p\!-\!|s'_k-s_k^{\ast}|^p\big|\notag\\
&\!\le\!\textstyle\sum_{k=1}^r\sum_{i=1}^p (^p_i)\Delta_{1,k}^{(i)}\left(\mbox{max}\{|c_{k0}|, |c_{k1}|\}\right)^{p-i}\!
\end{align}
\noindent When $p=1$,  Eqn (\ref{eqn:deltaup}) reduces to  $\Delta_u\le\sum_{k=1}^r\Delta_{1,k}$ in Part a). When $p=2$,  Eqn (\ref{eqn:deltaup}) becomes $\sum_{k=1}^r\!\left(\Delta_{1,k}^{(2)}+2\Delta_{1,k}\mbox{max}\{|c_{k0}|, |c_{k1}|\}\right)$, not as tight an upper bound as  Eqn (\ref{eqn:deltau2}). To see  this, we can show $2\Delta_{1,k}(c_{k1}-c_{k0})\le \Delta_{1,k}^{(2)}+2\Delta_{1,k}\mbox{max}\{|c_{k0}|, |c_{k1}|\}$ or  $2\Delta_{1,k}\mbox{max}\{|c_{k0}|, |c_{k1}|\}-2\Delta_{1,k}(c_{k1}-c_{k0})+\Delta_{1,k}^{(2)}\ge0$ holds for each $k$.
        When $c_{k0}c_{k1}\!\ge\!0$,  $c_{k1}-c_{k0}\!<\!\mbox{max}\{|c_{k0}|, |c_{k1}|\}$, $2\Delta_{1,k}(c_{k1}-c_{k0})\le 2\Delta_{1,k}\mbox{max}\{|c_{k0}|, |c_{k1}|\}<2\Delta_{1,k}\mbox{max}\{|c_\frac{}{}{k0}|, |c_{k1}|\}+\Delta_{1,k}^{(2)}$.
When $c_{k0}c_{k1}\le0$ and $\mbox{max}\{|c_{k1}|,c_{k0}|\}=c_{k1}$, $2\Delta_{1,k}\mbox{max}\{|c_{k0}|, |c_{k1}|\}-2\Delta_{1,k}(c_{k1}-c_{k0})+ \Delta_{1,k}^{(2)} =2\Delta_{1,k}c_{k1} -2\Delta_{1,k}(c_{k1}-c_{k0})+\Delta_{1,k}^{(2)}=2\Delta_{1,k}c_{k0}+\Delta_{1,k}^{(2)}$. \\ Since $\Delta_{1,k}^{(2)}\!=\!\!\maxxtex|s_k^2 -(s'_k)^2|\!=\!\maxxtex|s_k -s'_k|\cdot|s_k +s'_k| \ge \!\maxxtex\!|s_k -s'_k|\cdot|2c_{k0}|= 2\Delta_{1,k}|c_{k0}|$,  $\Delta_{1,k}^{(2)}-2\Delta_{1,k}|c_{k0}|=\Delta_{1,k}^{(2)}+2\Delta_{1,k}c_{k0}\ge0$.
    When $c_{k0}c_{k1}\le0$ and $\mbox{max}\{|c_{k1}|,c_{k0}|\}=|c_{k0}|$, $2\Delta_{1,k}\mbox{max}\{|c_{k0}|, |c_{k1}|\}+ \Delta_{1,k}^{(2)}-\Delta_{1,k}(c_{k1}-c_{k0})=2\Delta_{1,k}|c_{k0}|-2\Delta_{1,k}(c_{k1}-c_{k0})+\Delta_{1,k}^{(2)}=\Delta_{1,k}^{(2)}-2\Delta_{1,k}c_{k1}$. Since $\Delta_{1,k}^{(2)}=\!\!\maxxtex\!|s_k^2 -(s'_k)^2|\ge \!\!\!\!\maxxtex\!|s_k -s'_k|\cdot|2c_{k1}|= 2\Delta_{1,k}c_{k1}$,  $\Delta_{1,k}^{(2)}-2\Delta_{1,k}c_{k1}\ge0$.
All taken together,  $2\sum_{k=1}^r\Delta_{1,k}(c_{k1}-c_{k0})\le\sum_{k=1}^r\!\left(\Delta_{1,k}^{(2)}+2\Delta_{1,k}\mbox{max}\{|c_{k0}|, |c_{k1}|\}\right)$.
\end{proof}

\section{\normalsize{Proof of Lemma \ref{lem:lowerbound2}}}\label{app:lowerbound2}
\begin{proof}
When $r=1$ ($\s$ is a scalar), $\Delta_{p}\equiv\Delta$ for all $p\ge1$. To satisfy $(\epsilon,\delta)$-pDP, we set
\begin{align}\label{eqn:pDPr=1}
\!&\!\!\textstyle\Pr\!\left(|s^*\!-\!s|\!>\!\frac{\epsilon b^2\Delta^{-1}\!-\Delta}{2}\right)
\!=\!2\Phi\!\left(\!\frac{\Delta/2-\epsilon b^2(2\Delta)^{-1} }{b/\sqrt{2}}\!\right)\!\!\le\!\delta\\
&\Rightarrow\;\Delta b^{-1}\!\!-\!\epsilon b\Delta^{-1}\!\!\le\!\!\sqrt{2}\Phi^{-1}(\delta/2) \Rightarrow\;b\!\ge\!2^{-1/2}\epsilon^{-1}\Delta\sqrt{(\Phi^{-1}(\delta/2))^2+2\epsilon}-\Phi^{-1}(\delta/2).\notag
\end{align}
Together with the requirement $b^2\!-\!\epsilon^{-1}\Delta^2\!>\!0$, $b\!\ge\!\max\!\left\{\!\epsilon^{-1/2}\Delta, \left(\epsilon^{-1/2}\Delta\right)\!\frac{\sqrt{(\Phi^{-1}(\delta/2))^2+2\epsilon}-\Phi^{-1}(\delta/2)}{\sqrt{2\epsilon}}\!\right\}$.
Since $\delta\!<\!1$,  $\Phi^{-1}(\delta/2)\!<\!0$, $\sqrt{(\Phi^{-1}(\delta/2))^2\!+\!2\epsilon}\!- \! \Phi^{-1}(\delta/2) \ge\sqrt{2\epsilon}$, and thus \\
$b\!\ge\! \left(\epsilon^{-1/2}\Delta\right)\!\frac{\sqrt{(\Phi^{-1}(\delta/2))^2+2\epsilon}-\Phi^{-1}(\delta/2)}{\sqrt{2\epsilon}}$. When $r>1$, we leverage the proof in Appendix A  (page 265) in \cite{privacybook} and obtain
\begin{align*}
&\left|\log\!\!\left(\!\frac{\Pr(\s^{\ast} \in Q|\x)}{\Pr(\s^{\ast}\in Q|\x')}\!\right)\right|\!
=\!\left|\log\!\!\left(\!\frac{\exp\!\left(\!-\|\e\|_2^2/b^2\right) }{\exp\left(\!-\|\e+\mathbf{d}\|_2^2/b^2\right) }\!\right)\right|\\
=&\left|b^{-2}\left(\|\e\|_2^2-\|\e+\mathbf{d}\|_2^2\right) \right|
\le \left|b^{-2}\left(2\lambda\Delta_2+\Delta_2^2\right) \right|
\le b^{-2}\left( 2\Delta_2 |\lambda|+\Delta_2^2\right),
\end{align*}
where $\e=\s^*-\s(\x), \mathbf{d}=\s(\x)-\s(\x')$ defined in Eqn (\ref{eqn:second}), and $\lambda\sim N(0,b^2/2)$.  To satisfy $(\epsilon,\delta)$-pDP, we set
\begin{align*}
&\Pr(b^{-2}\left( 2\Delta_2 |\lambda|+\Delta_2^2\right)<\epsilon)= \Pr(\left(|\lambda|<(b^2\epsilon\Delta^{-1}_2-\Delta_2^2)/2\right)>1-\delta\\
\Rightarrow&\Pr\left(|\lambda|\!>\!(b^2\epsilon\Delta^{-1}_2\!\!-\!\Delta_2^2)/2\right)
\textstyle=2\Phi\!\left(\!\frac{\Delta_2-\epsilon b^2\Delta_2^{-1} }{\sqrt{2}b}\!\right)\!\!>\delta,
\end{align*}
which is the same as Eqn (\ref{eqn:pDPr=1}) for $r=1$. Similar to the case of $r=1$, we need $b^2\epsilon\Delta^{-1}_2-\Delta_2^2>0$, and the lower bound of $b$ for $r>1$ is
$b\!\ge\!\max\!\left\{\!\epsilon^{-1/2}\Delta_2, \left(\epsilon^{-1/2}\Delta_2\right)\!\frac{\sqrt{(\Phi^{-1}(\delta/2))^2+2\epsilon}-\Phi^{-1}(\delta/2)}{\sqrt{2\epsilon}}\!\right\}$. Since $\delta<1$,  $\Phi^{-1}(\delta/2)<0$, thus $b\!\ge\! \left(\epsilon^{-1/2}\Delta_2\right)\!\frac{\sqrt{(\Phi^{-1}(\delta/2))^2+2\epsilon}-\Phi^{-1}(\delta/2)}{\sqrt{2\epsilon}}$
 \end{proof}

\section{\normalsize{Proof of  Lemma \ref{lem:L1L2}}}\label{app:L1L2}
\begin{proof}
If $\sigma$ is set at the lower bound in Eqn (\ref{eqn:lowerbound2}), the ratio of the variance between the Gaussian distribution of the Gaussian mechanism of $(\epsilon,\delta)$-pDP  and the Laplace distribution of the  Laplace mechanism of $\epsilon$-DP is \begin{align}
&\textstyle\left(\!(2\epsilon)^{-1}\Delta_s\!\! \left(\!\sqrt{(\Phi^{-1}(\frac{\delta}{2}))^2+2\epsilon}\!-\!\Phi^{-1}(\frac{\delta}{2})\!\right)\!/(\!\sqrt{2}\epsilon^{-1}\Delta_s)\!\right)^{\!2}\notag\\
&\textstyle=\left(\!\sqrt{(\Phi^{-1}(\frac{\delta}{2}))^2+2\epsilon}-\Phi^{-1}(\frac{\delta}{2})\!\right)^2\!/8 =4^{-1}\!(\Phi^{-1}(\frac{\delta}{2}))^2\!+\!\epsilon\!-\! \Phi^{-1}(\frac{\delta}{2})\sqrt{(\Phi^{-1}(\frac{\delta}{2}))^2+2\epsilon}\label{eqn:varratio}
\end{align}
Since $\delta\!\in\![0,1]$,  $\delta/2\!\in\![0,0.5]$ and $\Phi^{-1}(\delta/2)\!\in\!(-\infty,0)$. Together with the fact $\epsilon\!>\!0$,  Eqn (\ref{eqn:varratio}) $\!>\!(\Phi^{-1}(\delta/2))^2/2$. Let $(\Phi^{-1}(\delta/2))^2/2\!>\!1$, then $\delta/2\!<\!\Phi(-\sqrt{2})$, leading to  $\delta<2\Phi(-\sqrt{2})\!\approx\!0.157$
\end{proof}
\end{appendix}

\bibliographystyle{ieeetran}
\input{arxivv5.bbl}

\end{document}